\documentclass[11pt]{amsart}

\usepackage[margin=1in]{geometry}
\usepackage{amsmath,amssymb,amsthm, comment, xcolor}
\usepackage{hyperref}
\usepackage{enumitem}

\theoremstyle{plain}
\newtheorem{theorem}{Theorem}[section]
\newtheorem{proposition}[theorem]{Proposition}
\newtheorem{corollary}[theorem]{Corollary}
\newtheorem{lemma}[theorem]{Lemma}
\newtheorem{conjecture}[theorem]{Conjecture}

\theoremstyle{definition}
\newtheorem{definition}[theorem]{Definition}
\newtheorem{example}[theorem]{Example}
\newtheorem{remark}[theorem]{Remark}
\newtheorem{question}[theorem]{Question}

\DeclareMathOperator{\wt}{wt}
\DeclareMathOperator{\Pf}{Pf}
\DeclareMathOperator{\GPf}{GPf}
\DeclareMathOperator{\Tr}{Tr}

\DeclareMathOperator{\Fun}{Fun}

\newcommand{\F}{\mathbb F}
\newcommand{\C}{\mathbb C}
\newcommand{\Z}{\mathbb Z}
\newcommand{\W}{W}
\newcommand{\RM}{\mathrm{RM}}

\newcommand{\lambdageo}{\lambda^{\operatorname{geo}}}

\subjclass[2020]{11E04, 11E08, 94B05, 94B27}

\keywords{quadratic form, Witt ring, powers of the fundamental ideal, iterated Laurent series field, gap theorems, Pfister numbers, symbol length, Reed--Muller codes, Kasami--Tokura classification, cancellation pattern}

\title[Quadratic Forms and Reed--Muller Codes]{Quadratic Forms and Reed--Muller Codes} 

\author{Neo Huang}
\address{Department of Mathematics\\
	University of British Columbia\\
	Vancouver, BC V6T 1Z2\\Canada}
    \email{nhuang04@student.ubc.ca}
\thanks{Neo Huang  was partially supported by an Undergraduate Student Research Award (USRA) from 
the Natural Sciences and Engineering Research Council of Canada.}

\author{Danny Ofek}   
\address{Department of Mathematics\\
	University of British Columbia\\
	Vancouver, BC V6T 1Z2\\Canada}
    \email{dannyofe@math.ubc.ca}
\thanks{Danny Ofek was supported by a postgraduate scholarship from the Natural Sciences and Engineering Research Council of Canada, reference number 601234.}

\author{Zinovy Reichstein}
\address{Department of Mathematics\\
	University of British Columbia\\
	Vancouver, BC V6T 1Z2\\Canada}
 \email{reichst@math.ubc.ca}
\thanks{{Zinovy Reichstein was partially supported by 
        an Individual Discovery Grant RGPIN-2023-03353 from the Natural Sciences and Engineering Research Council of Canada.}
}

\begin{document}

\begin{abstract}
The Witt ring of the field of iterated Laurent series $K_m = \C((t_1))\cdots((t_m))$
is isomorphic, as an $\F_2$-algebra, to the group algebra $\F_2[(\Z/2\Z)^m]$. Elements of this group algebra can, in turn, be identified with 
$\F_2$-valued functions on $\F_2^m$.
Under this identification, the $d$-th power $I^d(K_m)$ of the fundamental ideal corresponds
exactly to the binary Reed--Muller code $\RM(m,m-d)$, with scaled $d$-fold Pfister forms
corresponding to indicator functions of affine flats of dimension $d$. This
elementary observation, while built from ingredients that have 
been known since the 1960s and 1970s, does not appear to have been assembled into an
explicit dictionary between the two subjects. We record the dictionary and trace some of its consequences. In particular we show that
some of the classical results in the theory of Reed--Muller codes
can be deduced directly from the Arason--Pfister
Hauptsatz and the gap theorems in quadratic form theory. 
As an application of our dictionary, we prove new lower bounds on the maximal Pfister number of an anisotropic form of dimension $\leqslant n$ in $I^d(K_m)$.
\end{abstract}

\maketitle

\section{Introduction}

Two circles of ideas, developed independently and largely without reference to one
another, turn out to describe the same family of combinatorial objects. The first concerns
 the fundamental ideal $I(F)$ and its powers $I^d(F)$ in the Witt ring $\W(F)$ of a field $F$. These are central objects of study in the classical algebraic theory of quadratic forms over fields initiated by
Pfister in the 1960s with foundational contributions
by Milnor~\cite{Mil70} in 1970 and 
Arason--Pfister~\cite{AraPfi71} in 1971.
The second circle of ideas concerns a family of error-correcting codes $\RM(m,r)$ called  binary Reed--Muller codes. These codes were introduced in 1954 and studied
intensively by coding theorists from the 1960s onward, with foundational
contributions by Kasami, Lin, Peterson, Delsarte, Goethals, MacWilliams and
Tokura, among others. 

The bridge between the two is the field of iterated formal Laurent
series
\begin{equation} \label{e.Km}
  K_m \;=\; \C((t_1))(( t_2))\cdots((t_m)),
\end{equation}
where $\mathbb C$ denotes the field of complex numbers or, more generally, any quadratically closed field of
characteristic $\neq 2$. 
The Witt ring 
$\W(K_m)$ is isomorphic to the group algebra
$\F_2[V_m]$ of the elementary abelian group $V_m :=(\Z/2\Z)^m$ over the field $\F_2$ of two elements. This isomorphism, together with a description of the fundamental ideal as the
augmentation ideal, appears at least as early as the work of Knebusch, Rosenberg and
Ware in 1972 \cite{KneRosWar72} on the general problem of realizing Witt rings as
quotients of group rings. 
On the coding-theory side, it has been known since Berman's 1967 paper \cite{Ber67} that the
binary Reed--Muller code $\RM(m,m-d)$ is exactly the $d$-th power of the radical
(equivalently, the augmentation ideal) of $\F_2[V_m]$.

In Section~\ref{sec:dictionary}
we develop the resulting dictionary between quadratic form theory and the theory of Reed--Muller codes. Once assembled, the dictionary translates:
\begin{itemize}
  \item elements of $I^d(K_m)$ into codewords of $\RM(m,m-d)$,
  \item the dimension of the anisotropic part of a quadratic form into the Hamming
    weight of the corresponding codeword,
    \item $d$-fold Pfister forms into indicator functions of $d$-dimensional linear subspaces, 
  \item scaled $d$-fold Pfister forms into indicator functions of $d$-dimensional affine
    flats, and
  \item the $d$-Pfister number $\Pf_d(q)$ of a quadratic form $q$ with $[q] \in I^d(K_m)$ into 
  the minimum number of linear subspaces of dimension $d$ whose indicator functions sum to the
    corresponding codeword.
\end{itemize}

We emphasize that no single piece of the dictionary is difficult, and most of the raw ingredients have been available since the early 1970s. Our contribution is to make the connection between quadratic form theory and Reed--Muller codes explicit 
and trace the consequences in both directions. 

As we show in Section~\ref{sec:low_weights}, under this dictionary, the Arason--Pfister Hauptsatz \cite{AraPfi71} becomes exactly
the classical theorem that the minimum weight of $\RM(m,m-d)$ is $2^d$ and that the
minimum-weight codewords are precisely the indicator functions of single affine
$d$-flats. (This theorem is usually attributed, in the generality of $q$-ary codes, to
Delsarte, Goethals and MacWilliams \cite{DelGoeMac70}, building on the polynomial
formalism of Kasami, Lin and Peterson \cite{KasLinPet68}.) More strikingly, the Kasami--Tokura classification~\cite{KasTok70} 
of Reed--Muller codewords of low weight established the Gap Theorem for quadratic forms over $K_m$, 
three decades before the Gap Theorem was stated and proved over arbitrary fields 
by Hoffmann, Vishik and Karpenko~\cite{Hof98,Vis03,Kar04}. Moreover, while the Gap Theorem only predicts the possible 
dimensions of anisotropic forms $q$ in $I^d(F)$ in the range between $2^d$ and $2^{d+1} - 1$,
the Kasami--Tokura Theorem (reproduced as Theorem~\ref{thm:KT} below) translates under our dictionary to
a complete classification of such forms for $F = K_m$. 

We will use this classification in Section~\ref{sec:low_weight_decomp} 
to determine the exact value of the symbol lengths of such forms. In particular, we will establish the following.

\begin{theorem} \label{thm.length} {\rm(Theorem~\ref{thm:generalized_pfister_of_small})} Let $q$ be an anisotropic quadratic form over $K_m$ of dimension $2^{d+1} - 2^{d + 1 - \mu}$ such that $[q] \in I^d(K_m)$. Here $1 \leqslant \mu \leqslant d$. Then $\lambda_d(q) = \GPf_d(q) = \mu$ or $2$, depending on whether the word $w$ corresponding to $q$ belongs to family (1) or family (2) in the
Kasami--Tokura classification.
\end{theorem}

Here $\GPf_d(q)$ denotes the $d$th scaled Pfister number of $q$ and $\lambda_d(q)$ denotes the symbol length of the class of $q$ in $H^d(K_m, \Z/2\Z) = I^d(K_m)/I^{d+1}(K_m)$;
see Section~\ref{sec:pfister-numbers} for the precise definitions. 

We give several applications of this theorem in Section~\ref{sec:low_weight_decomp}. 
In particular, we give a counterexample to a question of Vishik (Question~\ref{q.vishik}) and state a modified version of this question as a conjecture
(Conjecture~\ref{q.vishik-modified}). Note that both Question~\ref{q.vishik} and Conjecture~\ref{q.vishik-modified} concern quadratic forms over an arbitrary field $F$. However, our contribution is based entirely on the dictionary constructed in Section~\ref{sec:dictionary}: our counterexample to Question~\ref{q.vishik} is over $F = K_m$, and our main piece of evidence for Conjecture~\ref{q.vishik-modified} is that it is true when $F = K_m$.

In the last two sections we use the dictionary to prove new lower bounds on Pfister numbers of forms defined over $K_m$. Our proofs are based on relating Pfister numbers to word-counts in $\RM(m,m-3)$. 
Let $F$ be a field of characteristic not two. Recall that the $d$-th Pfister number $\Pf_d(q)$ of a quadratic form $q\in I^d(F)$ is the minimal integer $s$ such that $q$ can be expressed as a sum $q = \sum_{i=1,\dots,s}\pm p_i$, where each $p_i$ is a $d$-fold Pfister form (see Section~\ref{sec:pfister-numbers}).
Thus, Pfister numbers measure the complexity of quadratic forms. Parimala, Suresh and Tignol \cite{ParSurTig09} described $\F_2[V]$ (for $V$ an $\F_2$-vector space) as ``a combinatorial analogue of the Witt ring'', identified Pfister elements, and used $\F_2[V]$ to study Pfister numbers of forms defined over the fields $K_m$. 
Their work raises the problem of computing the growth rate of 
$$\Pf_d(n) = \sup\{ \, \Pf_d(q) \, \mid \,  m \geqslant 1, \, q\in I^d(K_m),  \,  \dim(q) \leqslant n \, \} $$
as a function of $n$. Parimala, Suresh and Tignol computed $\Pf_1(n)$ and $\Pf_2(n)$ and showed that they were both linear in $n$. Raczek established the inequality $$\Pf_3(n) \leqslant n^2/8.$$ See \cite[Theorem 1.13]{Rac13}. It is known that $\Pf_3(q)$ can grow exponentially with the dimension of $q \in I^3(F)$ if $F$ is an arbitrary field \cite{BroReiVis10}.  Therefore Raczek's theorem shows that forms defined over $K_m$ are much simpler than forms over arbitrary fields.
Building on their work, Lorenz proved that $\Pf_d(n)$ grows at most like $n^{d-1}$ for all $d\geqslant  3$ \cite[Theorem 7.6]{Lor23}.  However, for forms over $K_m$ the best known lower bounds on $\Pf_d(n)$ are linear in $n$ for each $d \geqslant 2$ (in particular, $\displaystyle \Pf_3(n)\geqslant 2\lfloor \frac{n}{4}\rfloor -2$ for $n\geqslant 8$ by~\cite[Proposition 1.4]{Rac13}). 
In the last two sections of this paper we prove the first superlinear lower bound on $\Pf_d(n)$.

\begin{theorem}\label{thm:main}
For every $d \geqslant 3$ and every integer $n \geqslant 2^{d+4}$,
\[
  \Pf_d(n) >
  \frac{n^2}{9 \cdot 4^d \log_2(n)} .
\]
\end{theorem}
\noindent
See Section~\ref{sec:general_d} for a proof of a slightly more detailed version of Theorem~\ref{thm:main}.

\smallskip
Finally, a warning to the reader: this paper is not related to the line of research in coding theory in which the weight distribution of $\RM(m,2)$ is computed via the rank of the associated symplectic form (see, e.g.,~\cite{McE69}, \cite{SloBer70}, \cite[Chapter~15]{MacSlo77}). The quadratic form theory
used there concerns the classification of alternating forms over a finite field. Here our focus is on the structure of the Witt ring of the field $K_m$.

\subsection*{AI disclosure} In preparing this paper we made use of the AI assistants Gemini, ChatGPT and Claude. They helped us search the literature, proofread drafts, and 
construct some of the proofs.
All statements and proofs have been verified by the authors, who are solely responsible for the contents of this paper.

\section{Background}\label{sec:background}

\subsection{Quadratic forms and the Witt ring}

Let $F$ be a field of characteristic $\neq 2$. Write $\W(F)$ for the Witt ring of
$F$ and $I(F)$ for the fundamental ideal in $\W(F)$, i.e.,
the ideal of classes of even-dimensional forms. Recall that the powers
$I^d(F)$ are additively generated by the classes of $d$-fold \emph{Pfister forms}
\[
  \langle\!\langle a_1,\dots,a_d\rangle\!\rangle
  \;=\;
  \langle 1,-a_1\rangle\otimes\cdots\otimes\langle 1,-a_d\rangle,
  \qquad a_1,\dots,a_d\in F^\times.
\]
(We use the multiplicative convention $\langle\!\langle a\rangle\!\rangle = \langle
1,-a\rangle$; the sign is immaterial for our purposes since we will specialize to
fields in which $-1$ is a square.) For an introduction to these notions, we refer the reader to \cite{Lam05}.

The \emph{Arason--Pfister Hauptsatz} \cite{AraPfi71} states:

\begin{theorem}[Arason--Pfister, 1971]\label{thm:AP}
  Let $q$ be an anisotropic quadratic form over $F$ with $[q]\in I^d(F)$
  nonzero. Then $\dim(q) \geqslant 2^d$. If $\dim(q) = 2^d$, then $q$ is
  similar to a $d$-fold Pfister form.
\end{theorem}

In other words, $q$ is a non-zero scalar multiple of a $d$-fold Pfister form. We shall refer to such forms as scaled (or equivalently, generalized)
$d$-fold Pfister forms throughout this paper.

\subsection{Reed--Muller codes}

Fix $m \geqslant 0$. Identify a Boolean function $f:\F_2^m\to\F_2$ with the codeword
$(f(x))_{x\in\F_2^m}\in\F_2^{2^m}$. The $r$-th order binary \emph{Reed--Muller code}
$\RM(m,r)$, for $0\leqslant r\leqslant m$, is the $\F_2$-linear
span of the indicator functions of affine flats of dimension $d:=m-r$ (or codimension $r$). Equivalently,  $\RM(m, r)$ is the space of polynomial functions $\F_2^m \to \F_2$ of degree 
$\leqslant r$; see \cite[Chapter 13]{MacSlo77}.
For example, the polynomial function 
$f \colon \F_2^m \to \F_2$ of degree $2$
given by $f \colon (x_1, \ldots, x_m) \mapsto x_1 x_2$ is the indicator function of the $(m-2)$-dimensional affine flat given by $x_1 = x_2 = 1$. This function lies in $\RM(m, 2)$.
For an overview of the theory of Reed--Muller codes, we refer the reader to~\cite{ASY21}.

\begin{theorem}[Delsarte--Goethals--MacWilliams~\cite{DelGoeMac70}, 1970; Kasami--Lin--Peterson~\cite{KasLinPet68}, 1968]\label{thm:DGM}
  The minimum nonzero weight of $\RM(m,m-d)$ is $2^d$, and the minimum-weight
  codewords are exactly the indicator functions $\mathbf 1_D$ of affine flats $D$ of
  dimension $d$ (equivalently, of codimension $r=m-d$).
\end{theorem}

\begin{theorem}[Kasami--Tokura, 1970 \cite{KasTok70}]\label{thm:KT}
  Every nonzero codeword $w \in \RM(m,m-d)$ of weight strictly less than $2^{d+1}$
  has weight of the form $2^{d+1}-2^{d+1-\mu}$ for some integer $\mu\geqslant 1$ and
  is affine-equivalent to a Boolean function of one of the following two shapes:

\smallskip
(1) $w(x_1, \ldots, x_m) = x_1\cdots x_{r-2}\bigl(x_{r-1}x_r + x_{r+1}x_{r+2}+\cdots+
    x_{r+2\mu-3}x_{r+2\mu-2}\bigr)$, 
where $r + 2\mu - 2 \leqslant m$ and $\mu \geqslant 1$, or

\smallskip
(2) $w(x_1, \ldots, x_m) = x_1 \cdots x_{r - \mu}(x_{r-\mu + 1} \cdots x_r + x_{r+1} \cdots x_{r +\mu})$,  where $r + \mu \leqslant m$ and $r \geqslant \mu \geqslant 3$.

\smallskip \noindent 
Here $r := m-d$.
\end{theorem}

\subsection{Reed--Muller codes as radical powers}
\label{sect.berman}

Let $\F_2$ be the field of two elements, $V_m = (\Z/ 2 \Z)^m$ be an elementary abelian $2$-group of rank $m$,
and $\F_2[V_m]$ be the group algebra.
Denote a natural $\mathbb F_2$-basis element of $\F_2[V_m]$ corresponding to ${\bf v} \in V_m$ by
$X_{\bf v}$, so that $X_{\bf v} \cdot X_{\bf w} = X_{{\bf v} + {\bf w}}$. Let ${\bf e_1}, \ldots, {\bf e_m}$ be
the standard basis of $V_m$. Then
\begin{equation} \label{e.exterior}
\F_2[V_m] = \F_2[a_1, \ldots, a_m]/(a_1^2 - 1, \ldots, a_m^2 - 1) = \F_2[b_1, \ldots, b_m]/(b_1^2, \ldots, b_m^2) = \Lambda(V_m).
\end{equation}
Here $a_i := X_{\bf e_i}$, $b_i := 1 + a_i$, and $\Lambda(V_m)$ denotes the exterior algebra of $V_m$. (Note that $\Lambda(V_m)$ is commutative, because $V_m$ is a vector space over $\F_2$.)

We identify $\F_2[V_m]$ with the algebra $\Fun(\F_2^m, \, \F_2)$ of $\F_2$-valued functions on $V_m \cong \F_2^m$ as follows.
To an element $\sum_{{\bf v} \in V_m} c_{\bf v} X_{\bf v} \in \F_2[V_m]$ we associate
the function $\F_2^m \to \F_2$ given by ${\bf v} \mapsto c_{\bf v}$. Note that this identifies $\F_2[V_m]$ with
$\Fun(\F_2^m, \, \F_2)$ as $\F_2$-vector spaces. They are not isomorphic as $\F_2$-algebras. Indeed, $\F_2[V_m]$ has a large radical, whereas the radical of $\Fun(\F_2^m , \, \F_2)$ is trivial.

\begin{theorem}[Berman, 1967 \cite{Ber67}]\label{thm:Berman}
  Let $J = \Lambda^1(V_m) \oplus \cdots \oplus \Lambda^m(V_m) \subset \Lambda(V_m) = \F_2[V_m]$ 
  be the augmentation ideal.
   Under the identification of $\F_2[V_m]$ with the algebra
  of $\F_2$-valued functions on $V_m$ described above, the $d$-th power $J^d = \Lambda^d(V_m) \oplus  \cdots \oplus \Lambda^m(V_m)$ is exactly $\RM(m,m-d)$.
\end{theorem}

For reproofs and generalizations of Berman's theorem, see Charpin~\cite{Cha88} and Andriatahiny \cite{And16}.

\section{The dictionary}\label{sec:dictionary}

\subsection{The Witt ring of \texorpdfstring{$K_m$}{Km}}
Fix $m \geqslant 1$ and let $K_m = \C((t_1))((t_2))\cdots((t_m))$ be
an iterated Laurent series field, as in~\eqref{e.Km}. 
The square-class group of $K_m$ is
\begin{equation} \label{e.Vm}  V_m := K_m^\times / K_m^{\times 2} \;\cong\; (\Z/2 \Z)^m,
\end{equation}
generated by the classes of $t_1,\dots,t_m$. In other words, we identify $t^{\bf v} := t_1^{\epsilon_1} \ldots t_m^{\epsilon_m} 
\in K_m^{\times}/K_m^{\times 2}$ with ${\bf v} = (\epsilon_1, \ldots, \epsilon_m) \in (\Z/ 2 \Z)^m$.

\begin{theorem} 
\label{thm:springer}
  There is an isomorphism of $\F_2$-algebras
  \[
    \Theta: \W(K_m) \;\xrightarrow{\ \sim\ }\; \F_2[V_m] \cong \Lambda(V_m), \] 
  sending the class of the $1$-dimensional quadratic form $\langle t_i\rangle$
  to $X_{\bf e_i}$. 
\end{theorem}

The isomorphism between $\W(K_m)$ and $\F_2[V_m]$ is a consequence of Springer's theorem on the Witt ring of a complete
discretely valued field, applied recursively; see \cite[\S 2]{ParSurTig09}. The isomorphism between $\F_2[V_m]$ and $\Lambda(V_m)$ is explained in \eqref{e.exterior}.

\subsection{The dictionary} 
 Identify the group algebra $\F_2[V_m]$ with the space $\Fun(\F_2^m, \, \F_2)$ of Boolean functions $\F_2^m \to \F_2$, as in Section~\ref{sect.berman}. Let $\Theta$ be the isomorphism of Theorem~\ref{thm:springer}. 

\begin{theorem} \label{thm:dictionary}
 For any $0\leqslant d\leqslant m$:
  \begin{enumerate}[label=(\roman*)]
    \item $\Theta$ restricts to an isomorphism of $\F_2$-vector spaces
      $\displaystyle
        I^d(K_m) \;\xrightarrow{\ \sim\ }\; \RM(m,m-d)$.
    \item For ${\bf v_1},\ldots, {\bf v_d} \in V_m$ linearly independent and ${\bf a} \in V_m$,
    the class of the Pfister form $\langle\!\langle
      t^{\bf v_1},\dots,t^{\bf v_d}\rangle\!\rangle$ corresponds under $\Theta$ to the
      indicator function $\mathbf 1_D$ of the linear subspace
      $\displaystyle
        D = \mathrm{span}_{\F_2}({\bf v_1},\dots, {\bf v_d}) \subset \F_2^m$ and 
    the class of the scaled Pfister form $\langle t^{\bf a} \rangle\otimes\langle\!\langle
      t^{\bf v_1},\dots,t^{\bf v_d}\rangle\!\rangle$ to the
      indicator function $\mathbf 1_E$ of the affine flat
      $E = {\bf a} + D$.
    \item For an anisotropic form $q$ in $I^d(K_m)$, set
      $w =\Theta([q])$.
    Then
        $\dim (q) \;=\; \wt(w)$.
  \end{enumerate}
\end{theorem}

\begin{proof}
  (i) The isomorphism $\Theta$ of Theorem~\ref{thm:springer} carries the fundamental ideal $I(K_m)$ of $\W(K_m)$
   to the augmentation ideal $J\subset
  \F_2[V_m]$. Since $\Theta$ is a ring isomorphism, $\Theta(I^d(K_m)) = J^d$. By
  Theorem~\ref{thm:Berman}, $J^d = \RM(m,m-d)$ under the identification of
  $\F_2[V_m]$ with Boolean functions on $\F_2^m$.

  (ii) By definition,
  $\langle\!\langle t^{\bf v_1},\dots,t^{\bf v_d}\rangle\!\rangle = \bigotimes_i \langle 1,
  -t^{\bf v_i}\rangle$ has class, under $\Theta$, equal to $\prod_i (1+X_{\bf v_i}) \in
  \F_2[V_m]$. (Note that since $-1$ is a square in $K_m$, $-t^{\bf v_i}$ and $t^{\bf v_i}$ lie in the same square class.) Expanding,
  \[
    \prod_{i=1}^d (1+X_{\bf v_i}) = \sum_{S\subseteq\{1,\dots,d\}} X_{\sum_{i\in S} {\bf v_i}},
  \]
  which is exactly the indicator function of $D = \mathrm{Span}_{\F_2}({\bf v_1},\ldots, {\bf v_d})$ (a set of
  size $2^d$ since the $v_i$ are independent). Multiplying by the class of $\langle
  t^{\bf a}\rangle$, i.e.\ by $X_{\bf a}$ in $\F_2[V_m]$, translates the support by ${\bf a}$, giving the
  indicator function of $E = {\bf a} + D$.

  (iii) We claim that an $n$-dimensional quadratic form $q$ over $K_m$ is anisotropic if and only if 
  \begin{equation} \label{e.anisotropic}
  \text{$q \simeq \langle t^{\bf v_1},\dots, t^{\bf v_n} \rangle$ for some choice of $n$ distinct vectors  ${\bf v_1}, \ldots, {\bf v_n}$ in $V_m$.} 
  \end{equation}
  Once this claim is established, part (iii) readily follows. Indeed, since the support of $w = \Theta([q])$ is $\{ {\bf v_1}, 
  \ldots, {\bf v_n} \}$, we have $\dim(q) = n = \wt(w)$.
  
  \smallskip
  It remains to prove the claim. Since
  $K_m^{\times}/K_m^{\times 2}$ is generated by the classes of $t_1, \ldots, t_m$, every quadratic form $q$ is isomorphic to $\langle t^{\bf v_1},\dots, t^{\bf v_n} \rangle$ for some choice of ${\bf v_1}, \ldots, {\bf v_n}$ in $V_m$. If there is a repetition among
  these vectors, say ${\bf v_1} = {\bf v_2}$, then the $2$-dimensional subform $\langle t^{\bf v_1}, t^{\bf v_2} \rangle = \langle t^{\bf v_1} \rangle \otimes \langle 1, \, 1 \rangle$ is isotropic and hence so is $q$. (Here we use the fact that $-1$ is a square in $\C$ and hence, in $K_m$.) Hence, every anisotropic form $q$ is as in~\eqref{e.anisotropic}.

  For the converse, let us call the forms in~\eqref{e.anisotropic} ``candidates". {(Here we allow $\dim(q) = n$ to vary from $0$ to $2^m$.)} 
  As we have just seen, this pool of candidates contains every anisotropic form over $K_m$.  Our goal is to show that every one of these candidates is, in fact, anisotropic. We will do so by counting.  First we count the candidates. {Each candidate is uniquely determined by the choice of a subset $\{ {\bf v_1}, \ldots, {\bf v_n} \} \subset  V_m$.} Hence, 
  \[ \text{there are $2^{|V_m|}$ candidates.} \]
  On the other hand, by Theorem~\ref{thm:springer}, $|\W(K_m)| = 2^{|V_m|}$. Since each element of $\W(K_m)$ is represented by a unique anisotropic form, we conclude that 
  \[ \text{there are exactly $2^{|V_m|}$ anisotropic forms over $K_m$.} \]
   We conclude that {the candidates are distinct up to isomorphism, and each candidate is anisotropic.} This concludes the proof of the claim and thus Theorem~\ref{thm:dictionary}.
\end{proof}

\subsection{Milnor's conjecture and Reed--Muller codes}
\label{sect.milnor}

Milnor's conjecture asserts that for every field $F$ and every integer $d \geqslant 0$, 
there exists a surjective homomorphism of abelian groups $e_d \colon I^d(F) \to H^d(F, \mathbb Z/ 2 \mathbb Z)$ such that 
\begin{equation} \label{e.milnor}
\text{$\operatorname{Ker}(e_d) = I^{d+1}(F)$
and $e_d \colon 
\langle\! \langle a_1, \ldots, a_d \rangle\!\rangle 
\mapsto (a_1) \cdot \ldots \cdot (a_d) \in H^d(F, \Z/ 2 \Z)$}
\end{equation}
for any $d$-fold Pfister form $\langle\! \langle a_1, \ldots, a_d \rangle\! \rangle$ over $F$.

This conjecture was proved by Wadsworth~\cite{Wad83}
for strictly Henselian fields and by Voevodsky~\cite{Voe03} for arbitrary fields\footnote{More precisely, Voevodsky proved the K-theoretic version of Milnor's Conjecture in~\cite{Voe03}. The quadratic forms version we are 
discussing here was proved by Orlov, Vishik and Voevodsky in~\cite{OVV}.}.

In the case where $F = K_m = \mathbb C((t_1))\dots((t_m))$,
\[ \text{$\W(K_m) = \Lambda(V_m)$, and $J^d = \Lambda^d(V_m) \oplus \Lambda^{d+1}(V_m) \oplus \cdots \oplus \Lambda^m(V_m)$.} \]
Thus $\displaystyle J^d/J^{d+1} \simeq \Lambda^d(V_m)$, where $\Lambda^d(V_m)$ is in a natural bijection with $H^d(K_m, \Z/2\Z)$ given by $b_{i_1} \cdots b_{i_d} \mapsto (t_{i_1}) \cdots (t_{i_d})$. Here $b_1, \ldots, b_m$ are as in~\eqref{e.exterior}. In other words, the Milnor map is just the projection to the first summand
\[ e_d \colon J^d \to \Lambda^d(V_m) \simeq H^d(K_m, \Z/2\Z). \]

In the language of Boolean functions, an element 
of $\RM(m, m-d)$ is given by a polynomial function
$f \colon \F_2^m \to \F_2$ of degree $\leqslant 
r = m - d$. Any such function is uniquely represented as a sum of multilinear monomials of degree $\leqslant r$. The Milnor map $e_d$ takes 
a multilinear monomial of degree $\leqslant r-1$ to $0$, because any such monomial lies in $\RM(m, m-d -1)$, which corresponds to $I^{d+1}(K_m)$. On the other hand, we claim that for any multilinear monomial $x_{j_1} \ldots x_{j_r}$,
\begin{equation} \label{e.milnor-map} e_d \colon x_{j_1} \ldots x_{j_r} \mapsto (t_{i_1}) \cdot \ldots \cdot (t_{i_d}) \in H^d(K_m, \Z/ 2 \Z), 
\end{equation}
where $\{ j_1, \ldots, j_r \}$ and $\{i_1, \ldots, i_d \}$ are complementary subsets in $\{1, \ldots, m \}$. Indeed,
the monomial function $x_{j_1} \ldots x_{j_r} \colon \F_2^m \to \F_2$ is the indicator function 
of the $d$-dimensional flat $E$ given by $x_{j_1} = \dots = x_{j_r} = 1$. Let $D$ be the direction subspace for $E$, i.e., the $d$-dimensional $\F_2$-vector space $\operatorname{Span}_{\F_2}(e_{i_1}, \ldots, e_{i_d})$ given by $x_{j_1} = \dots = x_{j_r} = 0$. By Theorem~\ref{thm:dictionary}(ii), $E$ corresponds to the scaled Pfister form $\langle t^{\bf a} \rangle \otimes \phi$, where ${\bf a} := {\bf e_{j_1}} + \ldots + {\bf e_{j_r}} \in E$ and $\phi := \langle\!\langle t^{\bf e_{i_1}}, \dots, t^{\bf e_{i_d}} \rangle\!\rangle = \langle\!\langle t_{i_1}, \ldots, t_{i_d} \rangle\!\rangle$. Now observe that $\phi + \langle t^{\bf a} \rangle \otimes \phi = \langle 1, t^{\bf a} \rangle \otimes \phi \in I^{d+1}(K_m)$. By~\eqref{e.milnor},
$e_d(\langle t^{\bf a} \rangle \otimes \phi) = e_d(\phi) = (t_{i_1}) \cdots (t_{i_d})$, as claimed.

\begin{remark} \label{rem.hoffmann-becher} In $K$-theoretic language the Milnor isomorphism~\eqref{e.milnor-map} appears in \cite[(2.2) and Example 2.5]{BecHof04}.
\end{remark}

\begin{remark}  The quotient $\RM(m, r)/\RM(m, r-1)$
has received considerable attention in coding theoretic literature. Of particular interest is the minimal weight of a word in a given coset of $\RM(m, r-1)$ in $\RM(m, r)$. This is related to the covering radius of $\RM(m, r)$; for details and further references, see~\cite{Hou96, GL23}. In our interpretation, $\RM(m, r)/\RM(m, r-1)$ is identified with the cohomology group $H^d(K_m, \Z/ 2 \Z)$, where $d = m-r$. 
The minimal weight of a codeword representing a given element of $\RM(m, r)/\RM(m, r-1)$ is the 
same as the minimal dimension of an anisotropic quadratic form in $I^d(K_m)$ representing a given cohomology class in $H^d(K_m, \Z/ 2\Z)$. 
\end{remark}

\section{Low weights and the gap theorems}\label{sec:low_weights}
\subsection{Minimum distance}
\label{sec:min-distance}

As a first application, Theorem~\ref{thm:dictionary} relates the Arason--Pfister Hauptsatz to the Reed--Muller minimum-distance theorem.

\begin{corollary}\label{cor:AP-to-DGM}
   Theorem~\ref{thm:DGM} is equivalent to the case $F = K_m$ of Theorem~\ref{thm:AP} (Arason--Pfister).
\end{corollary}

\begin{proof}
  Let $0\neq w \in \RM(m,m-d)$, and let $\xi=\Theta^{-1}(w)\in I^d(K_m)$ be the
  corresponding nonzero class (Theorem~\ref{thm:dictionary}(i)). Let $q$ be an
  anisotropic representative of $\xi$. By Theorem~\ref{thm:dictionary}(iii),
  $\dim(q) = \wt(w)$. By the Hauptsatz, $\dim(q) \geqslant 2^d$, so $\wt(w) \geqslant 2^d$. Moreover,
  if $\dim(q) = \wt(w)=2^d$ then $q$ is a scaled
  $d$-fold Pfister form, which by
  Theorem~\ref{thm:dictionary}(ii) means $w$ is the indicator function of
  an affine $d$-flat. This proves Theorem~\ref{thm:DGM} for $\RM(m,m-d)$. The converse argument is identical, run backwards.
\end{proof}

\begin{remark} Note that the Delsarte--Goethals--MacWilliams theorem \cite{DelGoeMac70} is a statement about
generalized Reed--Muller codes over an arbitrary finite field $\F_l$.
Corollary~\ref{cor:AP-to-DGM} shows
that for $l = 2$ this theorem follows from the Arason--Pfister Hauptsatz.
Our argument does not recover the Delsarte--Goethals--MacWilliams
theorem for $l > 2$.
\end{remark}

We now turn to the second, and more striking, correspondence, concerning the next
threshold above the minimum weight.

\subsection{The gap theorems}\label{sect.gap}
On the quadratic-form side, the Gap Theorem  of Hoffmann, Vishik and Karpenko
\cite{Hof98,Vis03,Kar04} 
asserts that if $F$ is a field of characteristic $\neq 2$ and $q$ is an anisotropic form in $I^d(F)$ of dimension $\dim(q) < 2^{d+1}$, then $\dim(q) = 2^{d+1}-2^i$ for some $1\leqslant i\leqslant d$. Under the dictionary of
Theorem~\ref{thm:dictionary}, for $F = K_m$ this is the statement of
the Kasami--Tokura theorem (Theorem~\ref{thm:KT}) that the possible weights of
$\RM(m,m-d)$ below $2^{d + 1}$ are exactly $2^{d+1}-2^{d+1-\mu}$.
When $m \geqslant 2d$, the two families of
exceptional values match term for term. We are not aware of this numerical
coincidence having been remarked upon in either literature. 

As with Corollary~\ref{cor:AP-to-DGM}, we emphasize the asymmetry. On the one hand, the gap
theorem for a general field $F$ is considerably deeper than the gap statement 
in the Kasami--Tokura theorem, where $F = K_m$. On the other hand,
the Kasami--Tokura Theorem~\ref{thm:KT} gives a complete classification of words of weight $< 2^{d+1}$ in $\RM(m, m-d)$ 
or equivalently, of anisotropic forms $q$ of dimension $< 2^{d+1}$
over $K_m$ such that $[q] \in I^d(K_m)$.
This classification has no counterpart for general fields, 
except in the special case of forms of dimension $2^{d+1} - 2^{d-1}$, which we shall discuss in Section~\ref{sec:low_weight_decomp}; see
Conjecture~\ref{conj:lorenz}. {The classification
of words of weight $\leqslant 2^{d+1} + 2^{d-1}$
in $\RM(m, m-d)$ due to Kasami, Tokura and Azumi~\cite{KasTokAzu76} (a follow up to~\cite{KasTok70}) also has no counterpart for general fields.}

\section{Pfister numbers and the flat-decomposition problem}
\label{sec:pfister-numbers}

Let $F$ be a field. For a quadratic form $q$ over $F$ whose class lies in $I^d(F)$, we define

\begin{itemize}
\smallskip \item $\Pf_d(q)$, the Pfister number of $q$, as the minimal number $s$ such that $q$ is Witt equivalent to a signed sum $\sum_{i=1,\dots,s} \pm p_i$, where each $p_i$ is a $d$-fold Pfister form,

\smallskip \item
$\GPf_d(q)$, the generalized Pfister number of $q$, as the minimal number $s$ such that the Witt class of $q$ is a sum of $s$ scaled (or equivalently, generalized) $d$-fold Pfister forms in $I^d(F)$.

\smallskip \item
$\lambda_d(q)$, the symbol length of $q$, as the minimal number $s$ such that the class of $q$ is a sum of $s$ $d$-fold Pfister forms in $I^d(F)/I^{d+1}(F)$.
\end{itemize}

In the case where $F$ contains $\sqrt{-1}$ (and in particular, in the case where $F = K_m$, which is of particular interest to us throughout this paper), the signs can be absorbed into the forms $p_i$, and $\Pf_d(q)$ can be redefined as the minimal number $s$ such that $q$ is Witt equivalent to a direct sum of $s$ $d$-fold Pfister forms $p_1, \ldots, p_s$.

\smallskip
Recall from Section~\ref{sect.milnor} that $I^d(F)/I^{d+1}(F) \simeq H^d(F, \Z/ 2 \Z)$ via the Milnor isomorphism, so $\lambda_d(q)$ is the usual symbol length of the cohomology class in $H^d(F, \Z/ 2 \Z)$ represented by $q$. 

Note that we could have also defined the generalized symbol length $G \lambda_d(q)$ as the minimal number of scaled $d$-Pfister forms required to represent the class of $q$ in $H^d(F, \Z/2\Z)$, by analogy with $\GPf_d(q)$. The reason we did not do this is that a $d$-fold Pfister form
$\langle \!\langle a_1, \ldots, a_d \rangle \!\rangle$
represents the same class in $H^d(F, \Z/ 2\Z)$ as  any of its multiples $a \langle \!\langle a_1, \ldots, a_d \rangle \!\rangle$. Thus $G\lambda_d(q) = \lambda_d(q)$.

\begin{lemma} \label{lem.comparison}
Assume $[q] \in I^d(F)$. Then 
 $\displaystyle \lambda_d(q) \leqslant \GPf_d(q) \leqslant \Pf_d(q) \leqslant 2 \GPf_d(q)$.
\end{lemma}

\begin{proof} The first two inequalities are immediate from the definitions. The third 
inequality is proved 
in~\cite[Proposition 2.2]{Lor23}. 
\end{proof}

We now define analogous notions on the coding theoretic side. 
(This corresponds, in line with our dictionary, to setting $F = K_m$ in the above definitions.)
Given a word $w \in \RM(m, m-d)$, we denote by 

\begin{itemize}
\smallskip \item $l_d(w)$
the minimal number $s$ such that $w$ can be written as a sum of indicator functions of
$s$ $d$-dimensional subspaces in $\F_2^m$,

\smallskip \item $f_d(w)$,
the minimal number $s$ such that $w$ can be written as a sum of indicator functions of
$s$ $d$-dimensional affine flats in $\F_2^m$,

\smallskip \item
$s_d(w)$ the minimal number $s$ such that the class of $w$ in $\RM(m, m-d)/\RM(m, m-d-1)$
can be written as a sum of indicator functions of $s$ $d$-dimensional subspaces in $\F_2^m$.
\end{itemize}

Once again, we can replace $d$-dimensional subspaces by $d$-dimensional affine flats in the definition of $s_d(w)$. Indeed, modulo $\RM(m, m-d-1)$, the indicator function of a $d$-dimensional subspace equals the indicator function of any translate of this subspace. Also,
we have the following inequalities
\begin{equation} \label{e.comparison} s_d(w) \leqslant f_d(w) \leqslant l_d(w) \leqslant 2 f_d(w). \end{equation}
These inequalities can either be proved directly using combinatorial arguments or deduced from Lemma~\ref{lem.comparison} by substituting $K_m$ for $F$.

\begin{remark} \label{rem.independent-of-m}
For a given $w \in \RM(m, m-d)$, the numbers $l_d(w)$, $f_d(w)$ and $s_d(w)$ are ``independent of $m$" in the following sense. Suppose $w \colon V_m \to \F_2$ lies in $\RM(m, m-d)$, and let $w' \colon V_{m+1} \to \F_2$ be obtained from $w$ by extension by $0$, i.e., $w'(x_1, \ldots, x_{m+1}) = w(x_1, \ldots, x_m)(1+ x_{m+1})$. Then $w' \in \RM(m+1, m+1 - d)$
and
\[ l_d(w) = l_d(w'), \quad f_d(w) = f_d(w'), \quad \text{and}
\quad s_d(w) = s_d(w'). \]
Indeed, the inequalities $l_d(w) \geqslant l_d(w')$, $f_d(w) \geqslant f_d(w')$ and $s_d(w) \geqslant s_d(w')$ are obvious from the definitions: $d$-dimensional subspaces in $\F_2^m$ remain $d$-dimensional subspaces in $\F_2^{m+1}$ and similarly for flats.
This also shows that $w' \in \RM(m+1, m + 1 -d)$.

To prove the opposite inequalities, we use the push-forward map along the natural projection $\pi \colon \F_2^{m+1} \to \F_2^m$ taking a function $f \colon \F_2^{m+1} \to \F_2$ to the function $\pi_*(f) \colon \F_2^{m} \to \F_2$, where
$\pi_*(f)({\bf a})$ is the sum $f({\bf b_1}) + f({\bf b_2})$ over the fiber $\{ {\bf b_1}, \, {\bf b_2} \} = \pi^{-1}({\bf a})$. Then $\pi_*(w') = w$, since the restriction of $\pi$ to the supporting hyperplane $\{ x_{m+1} = 0 \}$ is an isomorphism between this hyperplane and $\F_2^m$. Note also that 
\[ \pi_*(1_D) = \begin{cases} 1_{\pi(D)}, \; \, \text{if $\dim(\pi(D)) = d $, and} \\
0, \; \, \text{if $\dim(\pi(D)) = d-1$}. \end{cases}
\]
Thus, if $w'$ can be written as a sum of indicator functions of $r$ linear subspaces (respectively, flats) $D_1, \ldots, D_r$, then
$w$ can be written as the sum of indicator functions of $\leqslant r$
linear subspaces (respectively flats), namely, those among $\pi(D_1), \ldots, \pi(D_r)$, which remain $d$-dimensional. This shows that $l_d(w) \leqslant l_d(w')$ and $f_d(w) \leqslant f_d(w')$. The argument for $s_d$ is similar.
\qed
\end{remark}

\smallskip
We now record the following immediate consequence of Theorem~\ref{thm:dictionary}.

\begin{corollary} \label{cor.l=Pf} Let $\xi \in I^d(K_m)$ and $w =\Theta(\xi)$ be the associated word in $\RM(m, m-d)$. Then $\Pf_d(\xi) = l_d(w)$ and $\GPf_d(\xi) = f_d(w)$ and $\lambda_d(\xi) = s_d(w)$. 
\qed
\end{corollary}
In the next section we will use Corollary~\ref{cor.l=Pf} in conjunction with the Kasami--Tokura theorem to compute the generalized Pfister numbers of low-dimensional forms in $I^d(K_m)$. 
 
Let $l_d(n)$ (respectively, $f_d(n)$ or $s_d(n)$) be the maximal value of
$l_d(w)$ (respectively, of $f_d(w)$ or of $s_d(w)$), as $w$ ranges over the words $w \in \RM(m, m-d)$ of weight $\leqslant n$. Similarly, let $\Pf_d(n)$ (respectively, $\GPf_d(n)$) be the maximal value of
$\Pf_d(q)$ (respectively, $\GPf_d(q)$), where $q$ is an anisotropic quadratic form in $I^d(K_m)$ of dimension $\leqslant n$. Here $d$ and $n$ are fixed but $m$ is allowed to range over the positive integers.
Corollary~\ref{cor.l=Pf} yields

\begin{corollary} \label{cor.l=Pf'} $\Pf_d(n) = l_d(n)$ and $\GPf_d(n) = f_d(n)$ for every $n \geqslant 1$. 
\qed
\end{corollary}

We can now translate some of the results on Pfister numbers into coding-theoretic language.

\begin{proposition} \label{prop.pfister-number} 
Let $n \geqslant 2$ be an even integer. Then 

\smallskip
(a) $l_2(n) = n-2$, and

\smallskip
(b) $\displaystyle f_2(n) = \frac{n}{2} - 1$. 

\smallskip
Moreover, for any integer $n \geqslant 1$, we have

\smallskip
(c) $\displaystyle l_3(n) \leqslant \frac{n^2}{8}$,

\smallskip
(d) $\displaystyle f_3(n) \leqslant \frac{n^2}{16}$.

\smallskip
(e) Now assume that $d \geqslant 2$ is an integer. Then there is a polynomial $p_d(n)$ of degree $d-1$ with rational coefficients such that $f_d(n) \leqslant p_d(n)$, for every positive integer $n$.
\end{proposition}

\begin{proof} 
(a) It is shown in~\cite{ParSurTig09} that $l_2(n) = n-2$. 

\smallskip
(b) By~\eqref{e.comparison},
$\displaystyle f_2(n) \geqslant \frac{1}{2} l_2(n) = \frac{n}{2} - 1$.
The opposite inequality follows from \cite[Exercise X.4]{Lam05}; cf.\ \cite[Proposition 2.3 and Example 2.4]{Lor23}. 

\smallskip
(c) and (d) are due to Raczek; see \cite[Theorem 1.13]{Rac13}. 
See also~\cite[Theorem 7.2]{Lor23} for a slightly stronger statement in (d).

\smallskip
(e) See~\cite[Theorem 7.6]{Lor23}. 
\end{proof} 

\begin{remark}
If we replace $K_m$ by an arbitrary field $F$, parts (c), (d) and (e) of Proposition~\ref{prop.pfister-number} fail dramatically. By \cite[Theorem 1.1]{BroReiVis10},
there is a sequence of fields $F_n$ and $n$-dimensional quadratic forms $q_n$ over $F_n$ such that $\Pf_3(q_n)$ increases exponentially with $n$.

If we fix $d \geqslant 4$,
it is not even known whether $\Pf_d(q)$ is bounded, as $F$ ranges over all fields and $q$ ranges over all forms of fixed dimension in $I^d(F)$.
\end{remark}

\begin{remark} Proposition~\ref{prop.pfister-number}(e) gives an upper bound on $s_d(n)$ as well, via $s_d(n) \leqslant f_d(n)$.
\end{remark} 

\begin{question} \label{q.pfister} 
What is the rate of growth of $s_d(n)$, $l_d(n)$ and $f_d(n)$ as $n \to \infty$?
\end{question}

\begin{remark} \label{rem.l-vs-f} In view of
\eqref{e.comparison}, $f_d(n) \leqslant l_d(n) \leqslant 2 f_d(n)$. In particular, for a fixed $d$,
$f_d(n)$ and $l_d(n)$ have the same growth rate as $n \to \infty$, up to a constant factor.
\end{remark}

\begin{remark} \label{rem.growth-rate}
It is not clear if the upper bound of Proposition~\ref{prop.pfister-number}(e) is optimal; the best lower bounds on $\Pf_d(n) = l_d(n)$ and $\GPf_d(n) =f_d(n)$ currently in the literature are linear in $n$; see \cite[Proposition 1.4]{Rac13} and \cite[Proposition 7.1]{Lor23}\footnote{These are stated for $d = 3$ only, but they can be leveraged to linear lower bounds on $\Pf_d(n)$ for every $d \geqslant 3$ by appealing to \cite[Proposition 3.3]{Lor23}. Our proof of Theorem~\ref{thm.general-d} will use \cite[Proposition 3.3]{Lor23} in a similar manner, to pass from $d -1$ to $d$.}. In Section~\ref{sec:lower_bd_l3},
we will show that $l_3(n)$ and $f_3(n)$ grow quadratically in $n$, up to a logarithmic factor. For $d \geqslant 4$, the gap between the best known upper and lower bounds on the asymptotic rate of growth of $l_d(n)$ is even wider (and similar for $f_d(n)$).
\end{remark}

\begin{remark} There is considerable literature on the symbol length in $H^d(F, \Z/ 2\Z)$ (and more generally in $H^d(F, \mu_p)$). The maximal value over all classes in $H^d(F, \Z/ 2 \Z)$ is an interesting  invariant of the field $F$. For symbol length in $H^d(K_m, \Z/ 2 \Z)$, see \cite{BecHof04, Kah05, HofLor25} and the references there. 

Usually the focus is on the maximal value of $s_d(w)$, as $w$ ranges over $\RM(m, m-d)$, though this maximal value is still unknown for $d \geqslant 3$.
(For $d = 2$, see \cite[Example 4.4(1)]{BecHof04}.) Here $m$ and $d$ are held fixed but $n = \wt(w)$ is allowed to vary. In this form, the problem can be rephrased as a question  about the minimal number of decomposable $d$-vectors needed to write any element of $\Lambda^d(V_m)$ as a sum. In other words, this is the $\F_2$-point version of the Grassmann secant variety problem. Our Question~\ref{q.pfister} is in a similar spirit, except that we are fixing $n$ and $d$ and allowing $m$ to vary.
\end{remark}

\begin{remark} We also note that Question~\ref{q.pfister} is distantly related to the line of research on the Schmidt rank and its variants, such as partition rank or Green--Tao rank. For an overview, see \cite{Pel24}. Roughly speaking, these results say that a lower bound on the bias of a word $w \in \RM(m, r)$ places an upper bound on the Schmidt rank of $w$. 
A lower bound on the bias is  
similar in spirit to requiring that $\wt(w) \leqslant n$; an upper bound on the Schmidt rank is similar - again, in spirit - to an upper bound on $f_d(w)$, where $d = m-r$. %
(More precisely, the Schmidt rank of the top-degree part of $w$ is at most $f_d(w)$.) 

The structure theory of biased polynomials \cite[Section 4]{Pel24} bounds Schmidt rank from above for polynomials of fixed degree in many variables; in our setting the degree $m - d$ is unbounded while $d$ is fixed, and these results, while similar in spirit, do not have a direct bearing on Question~\ref{q.pfister}.
\end{remark}
\section{Pfister numbers and symbol lengths in small dimensions}\label{sec:low_weight_decomp}

In Section~\ref{sect.gap}, we noted that the Kasami--Tokura theorem reproves the gap theorems for quadratic forms in $I^d(K_m)$ with dimension less than $2^{d+1}$. We will now use it to determine the generalized Pfister numbers and symbol lengths of such forms.

\begin{theorem}\label{thm:generalized_pfister_of_small}
     Let $q$ be an anisotropic quadratic form over $K_m$ with $[q]\in I^d(K_m)$ and assume $n=\dim q  < 2^{d+1}$. Recall that by the Kasami--Tokura Theorem~\ref{thm:KT} there exists $1\leqslant\mu\leqslant d$ such that $n = 2^{d+1} -2^{d+1-\mu}$. Moreover, the word $w = \Theta(q) \in \RM(m, m-d)$ associated to $q$ belongs to one of the two families, (1) and (2) given in Theorem~\ref{thm:KT}, up to affine equivalence.  Then
     \[
     \lambda_d(q) = \GPf_d(q) = \begin{cases} \mu, \; \; \text{if $\Theta(q)$ belongs to family (1),} \\
     2, \; \; \text{if $\Theta(q)$ belongs to family (2).} 
     \end{cases}
     \]
\end{theorem}

\begin{remark} \label{rem.mu<=2} If $w$ belongs to family (2), then $\mu \geqslant 3$. Thus the theorem tells us that if $\mu\leqslant 2$, then $\lambda_d(q) = \GPf_d(q) = \mu$.
\end{remark}

\begin{remark} \label{rem.one-family} Note that
Theorem~\ref{thm:generalized_pfister_of_small} implies that $w$ cannot belong to both families (1) and (2). If it did, then the theorem would tell us that $\mu = 2$. On the other hand, as we mentioned in Remark~\ref{rem.mu<=2}, if $w$ belongs to family (2), then $\mu\geqslant 3$, a contradiction. 
\end{remark}

\begin{proof}[Proof of Theorem~\ref{thm:generalized_pfister_of_small}]
The case where $n = 2^d$ and $\mu=1$ is covered by Theorem~\ref{thm:AP}. 
    Therefore we may assume that $n=2^{d+1}-2^{d+1-\mu}$ for some  $2\leqslant\mu\leqslant d$. By Corollary~\ref{cor.l=Pf} and Theorem~\ref{thm:dictionary}(iii), we have
    $$\wt(w) = n, \ f_d(w) =\GPf_d(q) \ \text{ and } \ s_d(w)=\lambda_d(q).$$
    Since affine changes of variables preserve $s_d$ and $f_d$, we may assume that $w$ is in one of the two normal forms described in cases (1) and (2) of Theorem~\ref{thm:KT}. Set $k=\mu$ in case (1) and $k=2$ in case (2). Since each monomial in the normal form is an affine $d$-flat indicator, we have $f_d(w)\leqslant k$. By \eqref{e.comparison}, it suffices to prove
    \begin{equation} 
    s_d(w) \geqslant k.
    \end{equation}
    Let $\pi:\F_2^m\to\F_2^{d+2}$ forget the first $r-2$ coordinates.
    Write elements of $\F_2^m$ as $(a,y)$, where
    $a\in\F_2^{r-2}$ and $y\in\F_2^{d+2}$.
    As in Remark~\ref{rem.independent-of-m}, we define the pushforward of a function $g:\F_2^m\to\F_2$ along $\pi$ by
\[
  \pi_*(g)(y):=\sum_{a\in\F_2^{r-2}}g(a,y).
\]
For any monomial $M=a_1^{i_1}\dots a_{r-2}^{i_{r-2}} y_1^{j_1}\dots y_{d+2}^{j_{d+2}}$ defined by integers $i_1,\dots,i_{r-2},j_1,\dots,j_{d+2}\geqslant 0$, we have
\begin{equation}\label{eq:formula_for_pi_*}
    \pi_*(M) = \begin{cases}
    y_1^{j_1}\dots y_{d+2}^{j_{d+2}} & \text{ if }i_t>0 \text{ for all }1\leqslant t\leqslant r-2 \\
    0 & \text{ otherwise.}  
\end{cases}
\end{equation}
In particular, $\pi_*$ sends polynomials of degree $r-1$ to affine-linear functions. Therefore
\begin{equation}\label{eq:pi_*_RM(m,r-1)}
\pi_*(\RM(m,r-1)) \subset \RM(d+2,1).
\end{equation}
Let $D\subset \F_2^m$ be an affine $d$-flat. We have
\begin{equation}\label{eq:pi_*_1_D}
  \pi_*(\mathbf1_D)=
  \begin{cases}
    \mathbf1_{\pi(D)},&\text{if }\dim\pi(D)=d,\\
    0,&\text{if }\dim\pi(D)<d
  \end{cases}
\end{equation}
because every nonempty fiber of $\pi|_D$ has
$2^{d-\dim\pi(D)}$ points. 

\smallskip
{\bf Claim:} If $\pi_*(w) = y_1 y_2 + \ldots + y_{2s-1} y_{2s}$, then $s_d(w) \geqslant s$.

\medskip
To prove the Claim, recall that by the definition of $s_d$, there exist affine $d$-flats $D_1,\dots,D_{s_d(w)}\subset \F_2^m$ such that
$$w \equiv  \mathbf{1}_{D_1} +\dots+ \mathbf{1}_{D_{s_d(w)}} \mod \RM(m,r-1).$$
Therefore \eqref{eq:pi_*_RM(m,r-1)} and \eqref{eq:pi_*_1_D} give
\begin{equation}\label{eq:pi_*_w_length}
\pi_*(w) = \mathbf{1}_{D'_1} +\dots +\mathbf{1}_{D'_t} +L,
\end{equation}
for some affine-linear function $L$, $t\leqslant s_d(w)$ and $d$-flats $D'_1,\dots,D'_t$. 

Let $U$ be the intersection of the direction spaces of the $d$-flats $D'_i\subset \F_2^{d+2}$.
Since $\operatorname{codim}(D'_i) =2$, we have $\operatorname{codim}(U)\leqslant 2t$. For every $u\in U$, the function $$\pi_*(w)(y+u)+\pi_*(w)(y) = L(y+u)+L(y)+ \sum_{i=1,\dots,t} \mathbf{1}_{D'_i}(y+u)+\mathbf{1}_{D'_i}(y)$$ is constant in $y\in \F_2^{d+2}$.
Since $\pi_*(w)$ is quadratic, expanding brackets gives
\[
  \pi_*(w)(y+u)+\pi_*(w)(y)
  =\pi_*(w)(u)+\sum_{i=1}^{s}
  \bigl(u_{2i-1}y_{2i}+u_{2i}y_{2i-1}\bigr).
\]
This forces $u_1=\cdots=u_{2s}=0$ for all $u\in U$. Hence
\[
  2s\leqslant\operatorname{codim}(U)
  \leqslant2t\leqslant2s_d(w).
\]
This completes the proof of the Claim.

\medskip
We will now consider two cases, according to the two normal forms in Theorem~\ref{thm:KT}.

\smallskip\noindent
\emph{Case (1):} $\displaystyle
     w(a,y)=a_1\cdots a_{r-2}
     \bigl(y_1y_2+\cdots+y_{2\mu-
 1}y_{2\mu}\bigr)$.
By \eqref{eq:formula_for_pi_*},
$\displaystyle
   \pi_*(w)=y_1y_2+\cdots+y_{2\mu-1}y_{2\mu}$.
The Claim now tells us that $s_d(w) \geqslant \mu$.

\smallskip \noindent
\emph{Case (2):} In this case, $r\geqslant \mu \geqslant3$ and $r+\mu \leqslant m$. Set $M_1(a) = a_1 \dots a_{r-\mu}$ and $M_2(a) = a_{r-\mu+1}\dots a_{r-2}$. By assumption, we have
$$w(a,y) = M_1(a)(M_2(a)y_1y_2 + y_3 y_4 \dots y_{\mu+2}).$$
Let $\tilde{w}\in \RM(m,m-d)$ denote the word obtained from $w$ by applying the linear substitution 
$$y_{i+2} \mapsto y_{i+2} +a_{r-\mu +i}\  \text{ for all }1\leqslant i\leqslant \mu-2.$$
Then $s_d(\tilde{w}) = s_d(w)$ and
$$\tilde{w}(a,y) =M_1(a)\big(M_2(a)y_1y_2 + y_{\mu+1}y_{\mu+2}\prod_{i=1,\dots,\mu-2}(y_{i+2} + a_{r-\mu+i}) \big).$$
By \eqref{eq:formula_for_pi_*}, applying $\pi_*$ extracts the $y$-coefficients of $M_1(a)M_2(a) = a_1\dots a_{r-2}$ from $\tilde{w}(a,y)$. Therefore
$$\pi_*(\tilde{w}(a,y)) = y_1y_2+ y_{\mu+1}y_{\mu+2}.$$
After relabeling the coordinates, so that $y_{\mu+1}, y_{\mu+2}$ become $y_3,y_4$, we can apply the Claim to conclude that
$s_d(w) = s_d(\tilde{w}) \geqslant 2$.
\end{proof}

Let $F$ be an arbitrary field with $\operatorname{char}(F) \neq 2$. The Arason--Pfister Hauptsatz shows that the case $n=2^d$ of Theorem~\ref{thm:generalized_pfister_of_small} holds with $K_m$ replaced by $F$. 
It is therefore natural to ask if the other cases of Theorem~\ref{thm:generalized_pfister_of_small} can be generalized to other fields as well. 

\subsection{The case $n= 2^{d+1}-2^{d-1}$}
The following conjecture is due to Hoffmann~\cite[Conjecture 2]{Hof98} for $d = 4$ and to Lorenz~\cite[Conjecture 1.5]{Lor25} for arbitrary $d \geqslant 2$.

\begin{conjecture}[Hoffmann, 1998; Lorenz, 2025] \label{conj:lorenz}
  Let $F$ be a field of characteristic $\neq 2$, let $d\geqslant 2$, and let $q$ be an
  anisotropic quadratic form over $F$ with $[q]\in I^d(F)$ and
  $\dim q = 2^{d+1}-2^{d-1}$. Then $[q]$ is
  a sum of two scaled $d$-fold Pfister forms in $W(F)$.
\end{conjecture}

An analogous statement for symmetric bilinear forms in characteristic $2$ was recently established by Scully \cite{Scu24}. Moreover, this conjecture was verified by Lorenz for fields of  (possibly transfinite) iterated Laurent series  $F = E((t_i))_{i\in I}$ over a linked field $E$;\footnote{The notion of a linked field is due to Elman and Lam~\cite{ElmLam72}. It means that any two $2$-fold Pfister forms $\pi_1$ and $\pi_2$ over $E$ can be written as $\pi_1 = \alpha \otimes \beta_1$ and $\pi_2 = \alpha \otimes \beta_2$ for some $1$-fold Pfister forms $\alpha, \beta_1, \beta_2$.} see \cite[Corollary 5.3]{Lor25}. 

We now observe that if we further specialize $F$ to $K_m$, then Conjecture~\ref{conj:lorenz} follows from
Theorem~\ref{thm:generalized_pfister_of_small}
and Remark~\ref{rem.mu<=2}. In fact, we don't need the full strength of Theorem~\ref{thm:generalized_pfister_of_small} here; the fact that $\GPf_d(q) \leqslant 2$
follows immediately from the $\mu=2$ case of the Kasami--Tokura Theorem~\ref{thm:KT}.
We also remark that in this special case our proof is very different from Lorenz's. His proof reduces, via his Theorem 5.1, to the base case of a linked field $E$, where the statement is vacuous: there are no anisotropic forms of dimension $2^{d+1}- 2^{d-1}$ in $I^d(E)$, where $E$ is a linked field; {see the proof of \cite[Corollary 5.3]{Lor25}}. The main content of his proof is in the  inductive step. For the fields $K_m$, we give a rather different direct non-inductive proof of Conjecture~\ref{conj:lorenz} by appealing to a combinatorial classification from 1970.

\subsection{The case $n = 2^{d+1}-2$ and Scharlau transfers}
\label{sect.scharlau}

Let $F$ be a field of characteristic different from $2$ and $E = F[\sqrt{\alpha}]$ be a quadratic field extension. Suppose $\phi$ is a quadratic form on a $n$-dimensional vector space over $E$. By the quadratic Scharlau transfer $\phi' = T_{E/F}(\phi)$ we shall mean the quadratic form $\phi' \colon V \times V \to F$ given by $v \to \operatorname{Tr}_{E/F}(\sqrt{\alpha} \phi(v))$. Here $\Tr_{E/F}$ denotes the trace function, and we view $V$ as an $F$-vector space of dimension $\dim_F(V) = 2 \dim_E(V)$. Note that $\alpha$ is uniquely determined by $E$, up to a square multiple from $F^{\times}$, so the Scharlau transfer, $\phi'$ is uniquely determined up to a scalar multiple.
 
This definition naturally extends to the case where $E/F$ is a quadratic \'etale algebra. The only quadratic \'etale algebra $E/F$ that is not a field extension is $E = F \oplus F$. Here we can also write $E = F[\sqrt{\alpha}]$, where
$\alpha = 1$ and $\sqrt{\alpha} = (1, -1) \in F \oplus F$. In this case a quadratic form $\phi$ over $E$ is a pair of quadratic forms $\phi_1, \phi_2$ of the same dimension over $F$, and the Scharlau transfer $T_{E/F}(\phi)$ is $\phi_1 \oplus (- \phi_2)$. By convention, a $d$-fold Pfister form over $E$ is a pair of $\phi = (\phi_1,\phi_2)$ where $\phi_1,\phi_2$ are $d$-fold Pfister forms over $F$.

The existence of $14$-dimensional forms $q\in I^3(F)$ such that $\GPf_3(q)=3$ over some field $F$ was established independently in \cite[Section 7]{izhboldin2000some} and \cite[Corollary 4.2]{hoffmann199814}. Therefore the case $n=14, d= 3$ of Theorem~\ref{thm:generalized_pfister_of_small} fails if we replace $K_m$ by the field $F$.  Nonetheless, the following theorem of Rost generalizes this case of Theorem~\ref{thm:generalized_pfister_of_small} to arbitrary fields. 

\begin{theorem}[\cite{Rost-14dim} {, see also \cite[Theorem 2.1(v)]{hoffmann199814}}]\label{thm:rost-transfer}
    Let $F$ be a field with $\operatorname{char}(F)\neq 2$. For every $14$-dimensional form $q \in I^3(F)$, there exists a quadratic \'etale $F$-algebra $E$ and a $3$-fold Pfister form $\phi$ such that $q$ is Witt equivalent to $T_{E/F}(\phi)$. 
\end{theorem}

Rost's theorem implies Theorem~\ref{thm:generalized_pfister_of_small} when $n=14,d=3$ because of the following lemma. 

\begin{lemma} \label{lem:generalized_pfister_of_transfer}
Let $E$ be a quadratic \'etale $K_m$-algebra for some $m\geqslant 1$. Assume $\phi$ is a scaled $d$-fold Pfister form over $E$. Then $q = T_{E/K_m}(\phi)$ is the sum of two scaled $d$-fold Pfister forms. 
\end{lemma}

\begin{proof} 
If $E = K_m\oplus K_m$ is the split quadratic algebra, then $T_{E/K_m}(\phi)$ is the difference of two scaled Pfister forms by definition. Hence we may assume $E = K_m[\sqrt{\alpha}]$ is a field. 
After multiplying $\alpha$ by squares, we may assume that $\alpha = t^{\bf v}$ for some $0 \neq {\bf v} = (\epsilon_1, \ldots, \epsilon_m) \in V_m$ without loss of generality. Let $i$ be the largest index such that $\epsilon_i = 1$. Then $\alpha = u t_i$ for some $u \in K_{i-1}^{\times}$ (where $K_0 = \C$), so $K_{i-1}((t_i)) = K_{i-1}((\alpha))$ and
\[ E = K_{i-1}((\sqrt{\alpha}))((t_{i+1})) \cdots ((t_m)) \]
is again an iterated Laurent series field over $\C$. In particular, $E^{\times}/E^{\times 2}$ is generated by the classes of $\sqrt{\alpha}$ and of $t_j$ for $j \neq i$. Thus every element of $E^{\times}$ lies in either $K_m^{\times} E^{\times 2}$ or $\sqrt{\alpha} K_m^{\times} E^{\times 2}$, and the product of any two elements of the latter set lies in the former.

Write $\phi = \langle c \rangle \otimes \langle\!\langle a_1, \ldots, a_d \rangle\!\rangle$ with $c, a_1, \ldots, a_d \in E^{\times}$. Since $-1$ is a square in $E$, we have $\langle\!\langle a, a' \rangle\!\rangle \simeq \langle 1, a, a', aa' \rangle \simeq \langle\!\langle aa', a' \rangle\!\rangle$ for any $a, a' \in E^{\times}$. Hence, after reordering $a_1, \ldots, a_d$ and replacing $a_j$ by $a_j a_d$ whenever neither $a_j$ nor $a_d$ lies in $K_m^{\times} E^{\times 2}$, we may assume that $a_1, \ldots, a_{d-1} \in K_m^{\times}$.

Now $\phi \simeq \psi_E \otimes \langle c, c a_d \rangle$, where $\psi = \langle\!\langle a_1, \ldots, a_{d-1} \rangle\!\rangle$ is a form over $K_m$. By Frobenius reciprocity~\cite[Section VII.1]{Lam05},
\[ T_{E/K_m}(\phi) \simeq \psi \otimes \bigl( T_{E/K_m}(\langle c \rangle) \oplus T_{E/K_m}(\langle c a_d \rangle) \bigr). \]
Each of the two transfers on the right is a binary form  over $K_m$ and every binary form is a scaled $1$-fold Pfister form. Therefore $T_{E/K_m}(\phi)$ is the sum of two scaled $d$-fold Pfister forms.
\end{proof}

Using Theorem~\ref{thm:generalized_pfister_of_small}, we will address the following question posed to us by A.~Vishik about a possible generalization of Rost's theorem (private communication).

\begin{question} \label{q.vishik} Let $F$ be a field of characteristic $\neq 2$, $d$ be a positive integer, and
$q$ be an $n$-dimensional non-degenerate quadratic form whose Witt class $[q]$ lies in $I^d(F)$.
Is it true that if $n < 2^{d+1}$, then $e_d(q) = e_d(T_{E/F}(\phi))$ in $H^d(F, \Z/2 \Z)$ for some 
quadratic \'etale $F$-algebra $E$ and some scaled $d$-fold Pfister form $\phi$ over $E$?
\end{question} 

To show that Question~\ref{q.vishik} has a negative answer, consider an anisotropic form $q\in I^d(K_m)$ such that $\Theta(q)$ comes from family (1) in the Kasami--Tokura classification (Theorem~\ref{thm:KT}) with $\mu = 3$, i.e., $\dim q =2^{d+1}- 2^{d-2}$. Such a form is available whenever $d\geqslant 4$. By Theorem~\ref{thm:generalized_pfister_of_small},
\begin{equation}\label{eq:vishik_eq}
\lambda_d(q) = 3.
\end{equation}
Meanwhile, if $e_d(q) = e_d(T_{E/K_m}(\phi))$ for some scaled $d$-fold Pfister form $\phi$ over an \'etale quadratic $K_m$-algebra $E$, then 
\begin{equation} \label{e.end-of-6.2} 3=\lambda_d(q) = \lambda_d(T_{E/K_m}(\phi)) \leqslant \GPf_d(T_{E/K_m}(\phi)) \leqslant 2, 
\end{equation}
a contradiction. (Here the first inequality uses Lemma~\ref{lem.comparison} and the second is Lemma~\ref{lem:generalized_pfister_of_transfer}.)
\qed

\subsection{Geometric symbol length}
\smallskip
We conclude this section with the following conjectural modification of Question~\ref{q.vishik} that remains open and is consistent with the results of this section. To state it, we use the notion of ``geometric symbol length'' due to Vishik. Let $q$ be a quadratic form over a field $F$. Assume that $\operatorname{char}(F)\neq 2$ and $[q] \in I^d(F)$. We define 
$$\lambdageo_d(q) = \min \Bigg\{ [L:F] \; \; \Bigg | \; \; \begin{matrix}
    e_d(q) = \operatorname{cor}_{L/F}(\alpha) \ \text{for some \'etale $F$-algebra $L$} \\ 
    \text{and some pure symbol
    $\alpha\in H^d(L,\Z/2\Z)$.}
\end{matrix} \   \Bigg\}$$
Here by a pure symbol we mean

\begin{itemize}
\smallskip \item the cup product of $d$ elements of
$H^1(L, \Z/2\Z)$ in the case where
$L/F$ is a finite separable field extension, and

\smallskip \item an $r$-tuple $(\alpha_1, \ldots, \alpha_r)$
in the case where $L = L_1 \oplus \ldots \oplus L_r$ is the direct
sum of $r$ finite separable field extensions $L_i/F$. Here each $\alpha_i \in H^d(L_i, \Z/ 2 \Z)$ is assumed to be 
a pure symbol over $L_i$.
\end{itemize}

Let $L/F$ be a finite separable field extension, and
$\alpha \in H^d(L, \Z/ 2\Z)$. We refer the reader to \cite[Section 99]{ekm} for the definition of the corestriction $\operatorname{cor}_{L/F}(\alpha) \in H^d(F, \Z/ 2 \Z)$. In particular, $\operatorname{cor}_{F/F}(\alpha) = \alpha$. If $L = L_1 \oplus \ldots \oplus L_r$ is a direct sum of $r$ finite separable field extensions $L_i/F$, we define the corestriction of $\displaystyle \alpha = (\alpha_1, \ldots, \alpha_r) \in H^d(L, \Z/ 2 \Z) = \bigoplus_{i = 1}^r H^d(L_i, \Z/ 2 \Z)$ as
\[ \operatorname{cor}_{L/F}(\alpha) = \sum_{i = 1}^r \operatorname{cor}_{L_i/F}(\alpha_i).\]
If $E$ is a quadratic \'etale $F$-algebra and
$\phi$ is a scaled $d$-fold Pfister form over $E$, then
\begin{equation} \label{e.cor}
\operatorname{cor}_{E/F}(e_d(\phi)) = e_d(T_{E/F}(\phi)).
\end{equation}
For a quadratic field extension $E/F$, this holds because the Scharlau transfer corresponds to the corestriction on mod $2$ Milnor $K$-theory by \cite[Lemma 40.1]{ekm}, and the norm residue homomorphism intertwines the corestriction on mod $2$ Milnor $K$-theory with the corestriction on Galois cohomology by \cite[Proposition 101.9]{ekm}.
For $E = F \oplus F$, formula~\eqref{e.cor} is immediate from the definitions. 

We also remark that in general the Scharlau transfer $T^{\,l}_{E/F}$ along a quadratic field extension  $E/F$
depends on the choice of a nonzero $F$-linear map $l \colon E \to F$. At the beginning of
Section~\ref{sect.scharlau} we chose this map to be $x \mapsto \Tr_{E/F}(\sqrt{\alpha}\,x)$.
However, if $\phi \in I^d(E)$, then the class
$e_d\bigl(T^{\,l}_{E/F}(\phi)\bigr) \in H^d(F,\Z/2\Z)$ does not depend on the
choice of $l$; see, e.g., \cite[Lemma 40.1]{ekm}. 

\smallskip
\begin{lemma} \label{lem.lambdageo-vs-lambda} Let $F$ be a field of characteristic different from $2$ and $q$ be a quadratic form over $F$ with $[q] \in I^d(F)$. Then $\lambdageo_d(q)\leqslant \lambda_d(q)$. 
\end{lemma}

\begin{proof} Set $r = \lambda_d(q)$ and $L = F^{\oplus r}$. 
 By definition, there exist pure symbols $(\alpha_1,\dots,\alpha_r) \in H^d(F,\Z/2\Z)^{\oplus r}$ such that $e_d(q) = \sum_{i=1,\dots,r}\alpha_i$. 
The corestriction of the pure symbol $\alpha = (\alpha_1,\dots,\alpha_r) \in H^d(L, \Z/ 2 \Z) =H^d(F,\Z/2\Z)^{\oplus r}$ 
is $\operatorname{cor}_{L/F}(\alpha) = \sum_{i=1,\dots,r}\alpha_i =  e_d(q)$. 
Therefore, $\lambdageo_d(q)\leqslant [F^{\oplus r}:F]= r = \lambda_d(q)$.
\end{proof}

\begin{conjecture} \label{q.vishik-modified} Let $F$ be a field of characteristic $\neq 2$, $d$ be a positive integer, and $q$ be an anisotropic quadratic form with $[q] \in I^d(F)$ such that  $\dim q<2^{d+1}$. Then $$\lambdageo_d(q) \leqslant 1+\lfloor d/2 \rfloor.$$
\end{conjecture}

\begin{remark} \label{rem1.vishik-modified}
Over $K_m$, this conjecture holds with the geometric length $\lambdageo_d(q)$ replaced by the ordinary symbol length $\lambda_d(q)$. Indeed,
in family~(1) of Theorem~\ref{thm:KT}, we have $\lambda_d(q)=\mu$
by Theorem~\ref{thm:generalized_pfister_of_small}. Since $r=m-d$,
the condition $r+2\mu-2\leqslant m$ gives $\lambda_d(q) = \mu\leqslant 1+\lfloor d/2\rfloor$.
In family~(2), $\lambda_d(q)=2$ and $d\geqslant\mu\geqslant 3$, so $\lambda_d(q) = 2 \leqslant 1 + \lfloor d/2 \rfloor$.  Since $\lambda_d^{\operatorname{geo}}(q)\leqslant \lambda_d(q)$ by Lemma~\ref{lem.lambdageo-vs-lambda}, we conclude that Conjecture~\ref{q.vishik-modified} holds for $F =K_m$.
\end{remark}

\begin{remark} \label{rem2.vishik-modified}
In the case where $q$ is our counterexample to Question~\ref{q.vishik} with $d = 4$, the inequality of Conjecture~\ref{q.vishik-modified} is, in fact, an equality. Recall that in this case $F = K_m$, $\mu = 3$, $\dim(q) = 2^{4+1} - 2^{4-2} = 28$, and $\lambda_4(q) = 3$ by Theorem~\ref{thm:generalized_pfister_of_small}. 
We claim that $\lambdageo_4(q) = 3$, which is exactly $1 + \lfloor d/2 \rfloor$.

By Lemma~\ref{lem.lambdageo-vs-lambda}, in order to establish this claim, we only need to show that $\lambdageo_4(q) > 2$. 
Assume the contrary. Then $e_4(q) = \operatorname{cor}_{L/K_m}(\alpha)$ for some
\'etale $K_m$-algebra $L$ of degree $[L: K_m] \leqslant 2$ and some pure symbol $\alpha \in H^4(L, \Z/ 2 \Z)$.
By definition of pure symbol, $\alpha = e_4(\phi)$ for a $4$-fold Pfister form $\phi$ over $L$. If $[L:K_m] = 1$, i.e., $L = K_m$, then $e_4(q) = \operatorname{cor}_{K_m/K_m}(e_4(\phi)) = e_4(\phi)$. Thus $\lambda_4(q) = 1$, a contradiction.
On the other hand, if $[L:K_m] = 2$, then $e_4(q) = \operatorname{cor}_{L/K_m}(\alpha) = \operatorname{cor}_{L/K_m}(e_4(\phi)) = e_4(T_{L/K_m}(\phi))$ by \eqref{e.cor}.
Arguing as in~\eqref{e.end-of-6.2}, we see that $3 = \lambda_4(q) \leqslant \GPf_4(T_{L/K_m}(\phi)) \leqslant 2$, a contradiction. \qed
\end{remark}

\begin{example} \label{ex2.vishik-modified}
For $d=3$, Conjecture~\ref{q.vishik-modified}  holds over every field $F$. Let $q$ be an anisotropic form of dimension $< 16$ with $[q] \in I^3(F)$. By the Gap Theorem (due to Pfister for $d = 3$), $\dim(q) = 8$, $12$ or $14$. If $\dim q\leqslant 12$, then results of Pfister give $\lambda_3(q) \leqslant \GPf_3(q)\leqslant 2$; see \cite[Theorem~2.1(i)--(iv)]{hoffmann199814}. If $\dim(q) = 14$,
Theorem~\ref{thm:rost-transfer} gives $[q]=[T_{E/F}(\phi)]$ for some scaled $3$-fold Pfister form $\phi$ over a quadratic \'etale $F$-algebra $E$.

If $E\cong F\oplus F$, then $q$ is the sum of two 
scaled $3$-fold Pfister forms, so $\lambdageo_3(q)\leqslant \lambda_3(q)\leqslant 2$. If $E/F$ is a quadratic field extension, then by~\eqref{e.cor}, we have 
$\displaystyle  e_3(q)=\operatorname{cor}_{E/F}(e_3(\phi))$. 
Here $e_3(\phi)$ is a single symbol; hence, $\lambdageo_3(q)\leqslant 2$. 
\end{example}

\begin{remark}
By \cite[Corollary 17.5]{izhboldin2000some}, there exists a field $F$ with $\operatorname{char}(F)\neq 2$ and a $14$-dimensional quadratic form $q$ with $[q] \in I^3(F)$ such that $\lambda_3(q)>2$. Thus one cannot hope to replace geometric symbol length $\lambdageo_d$ by the ordinary symbol length $\lambda_q$ in Conjecture~\ref{q.vishik-modified}.
\end{remark}

In summary, our counterexample to Question~\ref{q.vishik} shows that the constant $2 = [E:F]$ in Question~\ref{q.vishik} is too small. 
Conjecture \ref{q.vishik-modified} replaces $2$ with $1 + \lfloor d/2 \rfloor$, a natural candidate consistent with the evidence we have in the cases where 
(i) $F = K_m$ and $d$ is arbitrary and (ii) $d = 3$ and $F$ is arbitrary.

\section{A lower bound on \texorpdfstring{$l_3(n)$}{l3(n)}}\label{sec:lower_bd_l3}

In this section we will leverage our dictionary to prove the following lower bound on $\Pf_3(n)$. 

\begin{theorem}\label{thm:pfister_length_lower_bd}
    $\displaystyle \Pf_3(n) = l_3(n) \geqslant \frac{n^2}{256 \log_2(n)}$ for every even integer $n\geqslant 128$. 
\end{theorem}

The assumption that $n$ is even is harmless:
the lower bound on $\Pf_3(n)$ in Theorem~\ref{thm:pfister_length_lower_bd} extends to all integers $n \geqslant 128$ at the expense of replacing $256$ by a slightly larger constant.
Indeed, if $n$ is odd, then $\RM(m, m-3)$ has no words of weight $n$. Thus for $n$ odd and $\geqslant 129$, we have
\begin{align*}  l_3(n) = l_3(n-1) & \geqslant \frac{(n-1)^2}{256 \log_2(n-1)} >
\left( 1- \frac{1}{n} \right)^2  \cdot \frac{n^2}{256 \log_2(n)} \geqslant \left( 1- \frac{2}{n} \right)  \cdot \frac{n^2}{256 \log_2(n)} \\
 &  > \frac{63}{64} \cdot \frac{n^2}{256 \log_2(n)}  >  \frac{n^2}{260.1 \cdot \log_2(n)}.
\end{align*}

We also remark that the lower bound on $l_3(n)$ of Theorem~\ref{thm:pfister_length_lower_bd} is asymptotically close to the upper bound $\displaystyle l_3(n) \leqslant \frac{n^2}{8}$ of Proposition~\ref{prop.pfister-number}(c), though the two do not entirely match. Similarly, 
\[ \GPf_3(n) = f_3(n) \geqslant \frac{1}{2} l_3(n) \geqslant \frac{n^2}{512 \log_2(n)} \]
(see~\eqref{e.comparison}) is asymptotically close to the upper bound $\displaystyle f_3(n) \leqslant \frac{n^2}{16}$ of Proposition~\ref{prop.pfister-number}(d). Here $n$ is again an even integer $\geqslant 128$. In other words, Theorem~\ref{thm:pfister_length_lower_bd} shows that both $l_3(n)$ and $f_3(n)$ grow quadratically with $n$, up to a logarithmic factor.

Our proof of Theorem~\ref{thm:pfister_length_lower_bd} is entirely coding-theoretic; the dictionary is only used for the equality $\Pf_3(n) = l_3(n)$; see Corollary~\ref{cor.l=Pf'}.

For the remainder of this section, we fix an even integer $n \geqslant 128$ and set $l := l_3(n)$. 
Let $A_{n,m}$ denote the number of Hamming weight $n$ words in $\RM(m,m-3)$ for $m \geqslant n$. The proof of Theorem~\ref{thm:pfister_length_lower_bd} relies on comparing two estimates for $A_{n,m}$. The first estimate is based on the following elementary lemmas.

\begin{lemma}\label{lem:number_of_self_duals}
    There exist at least $2^{\binom{n/2}{2}-1}$ self-dual subspaces $U\subset \F_2^n$ that contain no vectors of Hamming weight $2$. 
\end{lemma}

\begin{proof}
    Let $T_{n}$  denote the number of self-dual subspaces of $\F_2^n$ and let $S_n$ be the number of self-dual subspaces of $\F_2^n$ which contain no words of weight $2$. By \cite[Section 2.2, Equation 5]{Rains-Sloane98}, 
    we have 
    \begin{equation}\label{eq:number_of_self_dual_subspaces}
        T_n = \prod_{i=1}^{n/2-1}(2^i+1).
    \end{equation}
    Note that we have
    \begin{equation}\label{eq:est_of_number_of_self_dual_subspaces}
        T_n \geqslant 2^{\sum_{i=1}^{n/2-1}i} = 2^{\binom{n/2}{2}}.
    \end{equation}
    For any $1\leqslant i< j \leqslant n$, there are exactly $T_{n-2}$ self-dual subspaces $V\subset \F_2^n$ that contain $w=e_i + e_j$.  Indeed, any such subspace decomposes as $V = \F_2 w\oplus V'$, where $V' \subset \F_2^{n-2}$ is a self-dual code and $\F_2^{n-2}$ denotes the $(n-2)$-dimensional subspace of $\F_2^n$ given by setting the $i$-th and the $j$-th coordinates equal to $0$. Moreover, $V'$ is uniquely determined by $w$ and $V$. As we vary over all pairs $1\leqslant i<j\leqslant n$, we see that there are at most $\binom{n}{2} T_{n-2}$ self-dual subspaces $V\subset \F_2^n$ that contain a weight $2$ word. 
    We compute
    $$S_n \geqslant T_n -\binom{n}{2}T_{n-2} \geqslant T_n(1-\frac{\binom{n}{2}}{2^{n/2-1}}) \geqslant \frac{1}{2}T_n \geqslant 2^{\binom{n/2}{2}-1}.$$ 
    Here the third inequality holds for all $n \geqslant 20$. 
    The last inequality is \eqref{eq:est_of_number_of_self_dual_subspaces}. The second inequality uses the inequality
    $$T_{n-2} = \frac{T_{n}}{2^{n/2-1}+1} \leqslant  \frac{T_{n}}{2^{n/2-1}},$$
    which follows from \eqref{eq:number_of_self_dual_subspaces}.
\end{proof}

\begin{lemma}\label{lem:elementary_count_of_spanning_tuples}
    Let $U$ be a $k$-dimensional $\F_2$-vector space. If $m > k$, then there are at least $2^{km-1}$ ways to choose $m$ vectors $u_1,\dots,u_m\in U$ that span $U$.
\end{lemma}

\begin{proof}
    Fix a basis for $U$. Choosing $m$ elements that span $U$ is equivalent to specifying a full rank $m$-by-$k$ matrix $X$. There are $2^m-2^{i-1}$ choices for the $i$-th column of $X$ if we choose them sequentially. Therefore overall, there are 
    $$\prod_{i=1}^{k}(2^m-2^{i-1}) = 2^{km} \prod_{i=1}^k (1-2^{i-1-m})$$
    ways to choose $u_1,\dots,u_m\in U$. Now the required lower bound follows from the elementary inequality
    $$\prod_{i=1}^k (1-2^{i-1-m}) \geqslant 1-2^{-m}\sum_{i=1}^k 2^{i-1}=1-\frac{2^k-1}{2^m} > 1-2^{k-m} \geqslant 1/2.$$
Here the last inequality relies on our assumption that $m > k$.
\end{proof}

\begin{lemma} \label{lem.lower-bound}
    For any integer $m\geqslant n$: 
    \begin{equation}\label{eq:estimate1_of_A_n_m}
    A_{n,m} \geqslant 2^{n m/2+n^2/16}.
    \end{equation}
\end{lemma}

\begin{proof}
    By \cite[Section II.C]{ASY21}, we have $\RM(m,m-3) = \RM(m,2)^{\perp}$. Thus, the indicator function $p = 1_{v_1,\dots,v_n}$ of $n$ distinct vectors is in $\RM(m,m-3)$ if and only if

    \smallskip
    \noindent
    (1) $\langle p,1\rangle =\sum_{i=1,\dots,n} 1= n$ is even (which we assume throughout this section),
    
    \smallskip
    \noindent
    (2) $\langle p,x_j\rangle  = \sum_{i=1\dots,n} v_i^j = 0 $ for all $1\leqslant j\leqslant m$, and

    \smallskip
    \noindent
    (3) $\langle p,x_jx_k\rangle = \sum_{i=1,\dots,n} v_i^j v_i^k= 0$ for all $1\leqslant j,k\leqslant m$.

\smallskip
Let $M$ be the $m$ by $n$ matrix whose columns are $v_1,\dots,v_n$. Conditions (2) and (3) are equivalent to $M M^t = 0$. Therefore codewords of weight $n$ in $\RM(m,m-3)$ correspond to matrices $M \in M_{m,n}(\F_2)$ whose columns are distinct and $M M^t=0$. 
In order to count such matrices, we shift attention to their row space. Let  $U_M \subset \F_2^n$ denote the row space of $M \in M_{m,n}(\F_2)$. The condition $M M^t = 0$ is equivalent to $U_M \subset U^{\perp}_M$, i.e. to $U_M$ being self-orthogonal. The columns of $M$ are distinct if and only if $U^{\perp}_M$ contains no weight $2$ vectors. By Lemma~\ref{lem:number_of_self_duals}, there are at least $2^{\binom{n/2}{2}-1}$ distinct self-dual subspaces $U\subset \F_2^n$ such that $U = U^{\perp}$ contains no weight $2$ vectors. From each such subspace $U\subset \F_2^n$, we can form at least $2^{mn/2 -1}$ matrices $M$ whose rows span $U$ by Lemma~\ref{lem:elementary_count_of_spanning_tuples} with $k = n/2$. Therefore there are at least
$$2^{\binom{n/2}{2}-1}\cdot 2^{nm/2 -1} =2^{\binom{n/2}{2}+nm/2 -2}$$
$m \times n$ matrices $M$ over $\F_2$ whose column vectors $v_1,\dots,v_n$ satisfy (1)-(3). Since the ordering of the column vectors was arbitrary, each indicator function $p=1_{v_1,\dots,v_n}$ corresponds to at most $n!$ of these matrices\footnote{In fact, each $p$ corresponds to precisely $n!$ such matrices but ``at most" suffices for the argument here.}. Therefore
$$A_{n,m} \geqslant \frac{2^{\binom{n/2}{2}+nm/2 -2}}{n!}.$$
To finish the proof, we apply the naive upper bound $\log_2(n!)\leqslant n \log_2(n)$ to obtain:
$$\log_2(A_{n,m}) \geqslant \binom{n/2}{2}+nm/2 -2 - n\log_2(n)=  n^2/8 +nm/2 - n/4 -2 -n\log_2(n).$$
An elementary computation shows that $n/4 +2 + n\log_2(n) \leqslant n^2/16$ for all $n\geqslant 128$. Therefore we conclude:
$$\log_2(A_{n,m}) \geqslant   n^2/8 +nm/2 - n^2/16 =  n^2/16 +nm/2.$$
\end{proof}

We will compare \eqref{eq:estimate1_of_A_n_m} to the following upper bound on $A_{n,m}$.
\begin{lemma}\label{lem:estimate2_of_A_n_m}
    For any integer $m\geqslant n$:
    \begin{equation}\label{eq:estimate2_of_A_n_m}
    A_{n,m} \leqslant 2^{n m/2}(8l)^{8l}.
    \end{equation}
\end{lemma}
The proof of \eqref{eq:estimate2_of_A_n_m} is based on counting words using cancellation patterns.

\begin{definition}
    Set $\Omega = \{1,\dots,l\}\times \F_2^3$. A \emph{cancellation pattern} is a pair $X=(S,\Pi)$, where $S \subset \Omega$ is a subset of cardinality $n$ and $\Pi$ is a partition of $\Omega\setminus S$ into pairs. We say that a word $p = (p_x)_{x\in \F_2^m} \in \RM(m,m-3)$ is \emph{compatible with $X$} if there exists a function $f:\Omega \to \F_2^m$, linear in the second coordinate such that

    \smallskip\noindent
    (1) The restriction $f_{|S}$ is injective.

    \smallskip\noindent
    (2) For any $\{x,y\}\in \Pi$, $f(x) = f(y)$.

    \smallskip\noindent
    (3) We have $p = 1_{f(S)}$, where $1_{f(S)}$ is the indicator function of $f(S) \subset \F_2^m$.

\smallskip\noindent A function $f$ satisfying these conditions will be called \emph{a parameterization function for $p$}.
\end{definition}

\begin{lemma}\label{lem:number_of_word_comptabile_with_pattern}
    Let $X = (S,\Pi)$ be a cancellation pattern. There are at most $2^{mn/2}$ codewords in $\RM(m,m-3)$ that are compatible with $X$.  
\end{lemma}
\begin{proof}
     Set $E = (\F_2^3)^{l}$. We view $\Omega$ as a subset of $E$ by identifying an element $(i,x)\in \Omega$ with the vector whose $j$-th coordinate is $x$ if $i=j$ and $0$ otherwise. We set $R_{X} = \langle x-y : \{x,y\} \in \Pi \rangle \subset E$ and $\overline{E}_{X} = E/ R_{X}$. Denote the quotient map $E\to \overline{E}_X$ by $v\mapsto \overline{v}$. The \emph{parameterization space} associated to $X$ is the linear subspace
    $$W_{X} =  \big\langle  \overline{u} \mid u\in S \big\rangle \subset \overline{E}_X .$$
    If $p$ is compatible with $X$, then there exists a parameterization function $f:\Omega \to \F_2^m$ for $p$. Extend $f$ linearly to a function $\tilde{F}: E\to \F_2^m$. Then $$\tilde{F}(x) - \tilde{F}(y) = f(x)-f(y) = 0$$ for all $\{x,y\} \in \Pi$. Therefore $\tilde{F}$ factors through $E\to \overline{E}_X$ and defines a map $F:\overline{E}_X\to \F_2^m$ such that
    $$F(\overline{x}) =f(x) \  \text{ for all } \ x\in \Omega.$$
    The word $p$ is uniquely determined by the restriction of $F$ to the parameterization space $W_X\subset \overline{E}_X$ because $p = 1_{f(S)}$ by definition. In particular, there are at most $$\lvert\operatorname{Hom}_{\F_2}(W_X,\F_2^m)\rvert = 2^{\dim_{\F_2}(W_X)m}$$ words $p\in \RM(m,m-3)$ compatible with $X$. It remains to prove
    \begin{equation}\label{eq:dim_of_param_space}
        \dim_{\F_2}(W_X) \leqslant n/2.
    \end{equation}
    In order to prove \eqref{eq:dim_of_param_space}, we use the embedding
    $$\operatorname{ev}_S:W_X^*\to \F_2^n, \ \ \lambda \mapsto (\lambda(\overline{u_1}),\dots,\lambda(\overline{u_n})).$$
    It suffices to prove that the image of $\operatorname{ev}_S$ is self-orthogonal with respect to the standard dot product~\footnote{The inequality \eqref{eq:dim_of_param_space} is a bit stronger than Raczek's~\cite[Proposition 1.2]{Rac13} but the proof proceeds along similar lines.}. 
    Let  $\mu,\lambda \in W_X^*$. We need to show
    \begin{equation}\label{eq:estimate2_A_n_m_final_goal}
        \sum_{i=1,\dots,n} \mu(\overline{u_i})\lambda(\overline{u_i}) = 0.
    \end{equation}
    Let $\tilde{\lambda},\tilde{\mu}\in \overline{E}_X^*$ be functionals that extend $\lambda$ and $\mu$ respectively. 
    For any $\{x,y\}\in \Pi$, we have $\overline{x} = \overline{y}$ in $\overline{E}_X$. Therefore
    $$\tilde{\lambda}(\overline{x})\tilde{\mu}(\overline{x}) + \tilde{\lambda}(\overline{y})\tilde{\mu}(\overline{y}) = 0.$$
    Since $\Pi$ partitions $\Omega\setminus S$, we get
    $$\sum_{j=1,\dots,l} \sum_{u\in \F_2^3}\tilde{\lambda}(\overline{j,u})\tilde{\mu}(\overline{j,u})=\sum_{x\in \Omega} \tilde{\lambda}(\overline{x})\tilde{\mu}(\overline{x})=\sum_{i=1,\dots,n} \mu(\overline{u_i})\lambda(\overline{u_i}).$$
    Therefore \eqref{eq:estimate2_A_n_m_final_goal} will follow once we prove
    \begin{equation}\label{eq:elementary_fact}\sum_{u\in \F_2^3}\tilde{\lambda}(\overline{j,u})\tilde{\mu}(\overline{j,u}) = 0
    \end{equation}
    for all $1\leqslant j\leqslant l$. Consider the linear functions $a(u) =\tilde{\lambda}(\overline{j,u})$ and $b(u)=\tilde{\mu}(\overline{j,u})$. 
    
    We claim that
    $\displaystyle \{ u \in \F_2^3 \, \mid \, a(u)=b(u)=1\}$ has an even number of elements.
    Indeed, this set is either empty or a coset of the
    vector subspace $\ker(a) \cap \ker(b)$ of $\F_2^3$ of dimension $i \geqslant 1$. Our set has 0 elements in the former case and $2^i$ elements in the latter. This proves the claim.

    Now
    $$\sum_{u\in \F_2^3}\tilde{\lambda}(\overline{j,u})\tilde{\mu}(\overline{j,u})=\sum_{u\in \F_2^3} a(u)b(u) \equiv |\{ u \mid a(u)=b(u)=1\}| \equiv 0 \mod 2.$$
    This establishes \eqref{eq:elementary_fact}, thus finishing the proof.
\end{proof}

\begin{proof}[Proof of Lemma~\ref{lem:estimate2_of_A_n_m}]
Since $l=l_3(n)$, any word $p\in \RM(m,m-3)$ of weight $n$ is the sum of $l'$ indicator functions of $3$-dimensional subspaces $V_1,\dots,V_{l'}$ of 
$\F_2^m$ for some $l' \leqslant l$.  
Choose a map $f:\Omega \to \F_2^m$ that restricts to an $\F_2$-linear isomorphism $\{i\}\times \F_2^3 \to V_i$
for $i \leqslant l'$ and to $0$ for $l' < i \leqslant l$.
Choosing one preimage of $f$ for each support point of $p$, we obtain a subset $S \subset \{ 1 , \ldots, l \} \times \F_2^3$ such that $f_{| \, S}$ is injective and
the images under $f$ of the elements of $\Omega \setminus S$ cancel in pairs. That is,
$p$ is compatible with some cancellation pattern $(S,\Pi)$ such that $f(x)=f(y)$ for all $\{x,y\}\in \Pi$. In other words, every word $p\in \RM(m,m-3)$ of weight $n$ is compatible with some cancellation pattern. By Lemma~\ref{lem:number_of_word_comptabile_with_pattern}, we obtain:
$$A_{n,m} \leqslant 2^{mn/2} C,$$
where $C$ is the number of cancellation patterns.
A cancellation pattern $X=(S,\Pi)$ is uniquely determined by the involution of $\Omega$ which exchanges all pairs appearing in $\Pi$ and fixes $S$. In particular, there are fewer patterns than the number of functions $\Omega \to \Omega$, which is $(8l)^{8l}$.   Therefore $C\leqslant (8l)^{8l}$ and \eqref{eq:estimate2_of_A_n_m} holds, as we wanted to show.
\end{proof}

We now have all the ingredients required to prove the main theorem of this section.
\begin{proof}[Proof of Theorem~\ref{thm:pfister_length_lower_bd}]
    Comparing \eqref{eq:estimate1_of_A_n_m} and \eqref{eq:estimate2_of_A_n_m} gives
    \begin{equation}\label{eq:proof_of_pfister_lower_bd_1}
    (8l)^{8l}\geqslant 2^{n^2/16}.
    \end{equation}
    
    If $l\geqslant n^2/8$, then there is nothing to prove. Therefore we may thus assume without loss of generality that $l\leqslant n^2/8$. Now \eqref{eq:proof_of_pfister_lower_bd_1} yields
    $$n^2 \leqslant 128l\log_2(8l)\leqslant 128l\log_2(n^2) = 256l\log_2(n), $$
    as desired.
\end{proof}

\begin{remark} \label{rem.exhibit} Previous (linear)
lower bounds on $\Pf_3(n)$ and $\GPf_3(n)$ were proved constructively, by exhibiting a specific $n$-dimensional quadratic form $q$ and bounding $\Pf_3(q)$ from below; see~\cite[Proposition 1.4]{Rac13} and \cite[Proposition 7.1]{Lor23}. Our proof of Theorem~\ref{thm:pfister_length_lower_bd} is non-constructive: we do not prove that $\displaystyle \Pf_3(q) \geqslant \frac{n^2}{256 \log_2(n)}$ for any specific quadratic form $q$ of dimension $\leqslant n$ whose class lies in $I^3(K_m)$. We only show that such a $q$ exists.  
\end{remark}

\begin{remark}
It is tempting to  try to modify the proof of Theorem~\ref{thm:pfister_length_lower_bd} to show that
$\Pf_d(n)$ is bounded from below by a polynomial of degree $d-1$, up to a logarithmic factor; see Remark~\ref{rem.growth-rate}.
The main obstacle is in extending Lemma~\ref{lem.lower-bound}. For general $d \geqslant 3$ the row space $U_M$ becomes $(d-1)$-fold self-orthogonal. For $d \geqslant 4$ such codes are much more constrained than self-orthogonal codes: there are fewer of them,  their structure is poorly understood, and the condition that a $(d-1)$-fold self-orthogonal code does not contain a vector of weight $2$ becomes more stringent. 
Informally speaking, $(d-1)$-fold self-orthogonal codes behave in a similar manner to linear subspaces on hypersurfaces of degree $d-1$ in $\mathbb P^{n-1}$ over an algebraically closed field of characteristic $0$. Linear subspaces on a quadric are plentiful, and each of them lies in one of maximal dimension
(i.e., projective dimension $\displaystyle \frac{n}{2} - 1$). Linear subspaces on higher degree hypersurfaces tend to be of smaller dimension and more sparsely distributed.
\end{remark}

However, in the next section we will see that Theorem~\ref{thm:pfister_length_lower_bd} can be leveraged to obtain a quadratic lower bound (up to a logarithmic factor) on $\Pf_d(n) = l_d(n)$ for every
$d \geqslant 3$.

\section{A lower bound on 
\texorpdfstring{$l_d(n)$}{ld(n)}}\label{sec:general_d}
\begin{theorem} \label{thm.general-d}
For every $d \geqslant 3$ and every integer $n \geqslant 2^{d+4}$,
\[
  \Pf_d(n) \;=\; l_d(n) \;\geqslant\; \GPf_d(n) \;=\; f_d(n) >
  \frac{n^2}{9 \cdot 4^d \log_2(n)} .
\]
\end{theorem}

The equalities $\Pf_d(n) = l_d(n)$ and $\GPf_d(n) = f_d(n)$ are taken from Corollary~\ref{cor.l=Pf'}. We will also make use of the inequalities $\displaystyle l_d(n) \geqslant f_d(n) \geqslant \frac{1}{2} l_d(n)$; see~\eqref{e.comparison}.

We begin with the following lemmas.

\begin{lemma} \label{lem.lorenz-prop3.3}
For every $d \geqslant 2$ and every even $n \geqslant 2$, we have $f_d(n) \geqslant f_{d-1}(n/2)$.
\end{lemma}

\begin{proof} Choose a quadratic form $q$ of dimension $\leqslant n/2$ with $[q] \in I^{d-1}(K_m)$ 
such that 
\[ \GPf_{d-1}(q) = \GPf_{d-1}(n/2) = f_{d-1}(n/2). \]   
Then $q'=\langle\!\langle t_{m+1}\rangle\!\rangle\otimes q$ is a form of dimension $\dim(q') = 2 \dim(q) \leqslant n$ with
$[q'] \in I^d(K_{m+1})$\footnote{Recall that $K_{m+1}$ is defined as $K_m((t_{m+1}))$.}. 
By \cite[Proposition 3.3]{Lor23}, $f_d(n) = \GPf_d(n) \geqslant \GPf_d(q') = \GPf_{d-1}(q)= \GPf_{d-1}(n/2) = f_{d-1}(n/2)$. 
\end{proof}

\begin{lemma} \label{lem.floor} 
For every $d \geqslant 4$ and every integer $N \geqslant 1$,
\[
  f_d(N) \;\geqslant\; f_{d-1}\bigl(\lfloor N/2\rfloor\bigr) .
\]
\end{lemma}

\begin{proof}
Since $2\lfloor N/2\rfloor \leqslant N$ and $f_d$ is non-decreasing,
$f_d(N) \geqslant f_d\bigl(2\lfloor N/2\rfloor\bigr)$. The integer
$n = 2\lfloor N/2\rfloor$ is even, so Lemma~\ref{lem.lorenz-prop3.3} applies to it and gives
$f_d\bigl(2\lfloor N/2\rfloor\bigr) \geqslant
f_{d-1}\bigl(\lfloor N/2\rfloor\bigr)$.
\end{proof}

\begin{proof}[Proof of Theorem~\ref{thm.general-d}]
Iterating the inequality of Lemma~\ref{lem.floor} $d-3$ times, and using the identity
$\bigl\lfloor\lfloor N/a\rfloor/b\bigr\rfloor = \lfloor N/ab\rfloor$, we arrive at
$f_d(n) \geqslant f_3\bigl(\lfloor n/2^{\,d-3}\rfloor\bigr)$. After one further
monotone step, $\lfloor n/2^{\,d-3}\rfloor \geqslant
2\lfloor n/2^{\,d-2}\rfloor$, we land on an even integer:
\begin{equation} \label{e.chain}
  f_d(n) \;\geqslant\; f_3\bigl(2\lfloor n/2^{\,d-2}\rfloor\bigr).
\end{equation}
If $n \geqslant 2^{d+4}$, then 
\begin{equation} \label{e.domain} n/2^{d-2} \geqslant 64.
\end{equation}
Hence,
$\displaystyle 2\lfloor n/2^{\,d-2}\rfloor \geqslant 128$, and 
we can apply Theorem~\ref{thm:pfister_length_lower_bd}:
\begin{equation} \label{e.chain2}
f_d(n) \geqslant f_3 \big( 2\lfloor n/2^{\,d-2}\rfloor \big) \geqslant \frac{1}{2} l_3 \big( 2\lfloor n/2^{\,d-2}\rfloor \big) \geqslant \frac{2 \big(\lfloor n/2^{\,d-2}\rfloor\big)^2}{256 \log_2 \big( 2\lfloor n/2^{\,d-2}\rfloor\big)} \geqslant  \frac{\big(\lfloor n/2^{\,d-2}\rfloor\big)^2}{2^7  \log_2 \big( 2\lfloor n/2^{\,d-2}\rfloor\big)}. 
\end{equation}
Using the inequality $\displaystyle 2\lfloor n/ 2^{d-2} \rfloor \leqslant n$, we can replace $\displaystyle \log_2 \big(2\lfloor n/2^{\,d-2}\rfloor\big)$ by $\log_2(n)$. Moreover, remembering \eqref{e.domain}, we see that
\[ \big( \lfloor \frac{n}{2^{d-2}} \rfloor \big)^2 \geqslant (\frac{n}{2^{d-2}} - 1)^2 \geqslant  \frac{n^2}{2^{2d-4}}(1 - \frac{2^{d-2}}{n})^2 \geqslant \frac{n^2}{2^{2d-4}} \cdot \big(1 - 2 \cdot \frac{2^{d-2}}{n} \big) \geqslant  \frac{n^2}{2^{2d-4}} \cdot \big( 1-\frac{1}{32} \big). \]
Substituting this into~\eqref{e.chain2}, we obtain
\[ f_d(n) \geqslant \frac{31}{32} \cdot \frac{n^2}{2^{2d+3}  \log_2(n)} >  \frac{n^2}{9 \cdot 4^d  \log_2(n)}, \]
where the last inequality follows from $\displaystyle \frac{31}{32} > \frac{8}{9}$.
\end{proof}

\section*{Acknowledgments} We are grateful to Detlev Hoffmann, David Leep and Alexander Vishik for helpful comments.


\providecommand{\bysame}{\leavevmode\hbox to3em{\hrulefill}\thinspace}
\providecommand{\MR}{\relax\ifhmode\unskip\space\fi MR }
\providecommand{\MRhref}[2]{%
  \href{http://www.ams.org/mathscinet-getitem?mr=#1}{#2}
}
\providecommand{\href}[2]{#2}


\begin{thebibliography}{KRW72}

\bibitem[And22]{And16}
Harinaivo Andriatahiny, \emph{On algebraic structure of the {R}eed-{M}uller
  codes}, Art Discrete Appl. Math. \textbf{5} (2022), no.~2, Paper No. 2.07, 9.
  \MR{4448657}

\bibitem[AP71]{AraPfi71}
J{\'o}n~Kristinn Arason and Albrecht Pfister, \emph{Beweis des {K}rullschen
  {D}urchschnittsatzes f{\"u}r den {W}ittring}, Inventiones Mathematicae
  \textbf{12} (1971), 173--176.

\bibitem[ASY21]{ASY21}
Emmanuel Abbe, Amir Shpilka, and Min Ye, \emph{Reed-{M}uller codes: theory and
  algorithms}, IEEE Trans. Inform. Theory \textbf{67} (2021), no.~6, part 1,
  3251--3277. \MR{4289318}

\bibitem[Ber67]{Ber67}
S.~D. Berman, \emph{On the theory of group codes}, Kibernetika \textbf{\rm
  (Kiev)} (1967), no.~1, 31--39. \MR{270816}

\bibitem[BH04]{BecHof04}
Karim~Johannes Becher and Detlev~W. Hoffmann, \emph{Symbol lengths in {M}ilnor
  {$K$}-theory}, Homology Homotopy Appl. \textbf{6} (2004), no.~1, 17--31.

\bibitem[BRV10]{BroReiVis10}
Patrick Brosnan, Zinovy Reichstein, and Angelo Vistoli, \emph{Essential
  dimension, spinor groups, and quadratic forms}, Annals of Mathematics
  \textbf{171} (2010), 533--544.

\bibitem[Cha88]{Cha88}
Pascale Charpin, \emph{Une g\'{e}n\'{e}ralisation de la construction de
  {B}erman des codes de {R}eed et {M}uller {$p$}-aires}, Comm. Algebra
  \textbf{16} (1988), no.~11, 2231--2246. \MR{962313}

\bibitem[DGM70]{DelGoeMac70}
P.~Delsarte, J.-M. Goethals, and F.~J. MacWilliams, \emph{On generalized
  {R}eed--{M}uller codes and their relatives}, Information and Control
  \textbf{16} (1970), 403--442.

\bibitem[EKM08]{ekm}
Richard Elman, Nikita Karpenko, and Alexander Merkurjev, \emph{The algebraic
  and geometric theory of quadratic forms}, American Mathematical Society
  Colloquium Publications, vol.~56, American Mathematical Society, Providence,
  RI, 2008. \MR{2427530}

\bibitem[EL72]{ElmLam72}
Richard Elman and T.~Y. Lam, \emph{Pfister forms and {$K$}-theory of fields},
  J. Algebra \textbf{23} (1972), 181--213. \MR{302739}

\bibitem[GL23]{GL23}
Val\'{e}rie Gillot and Philippe Langevin, \emph{Classification of some cosets
  of the {R}eed-{M}uller code}, Cryptogr. Commun. \textbf{15} (2023), no.~6,
  1129--1137. \MR{4663601}

\bibitem[HL25]{HofLor25}
Detlev~W. Hoffmann and Nico Lorenz, \emph{On the number of {P}fister forms and
  the symbol length of fields with finite square class number}, Homology
  Homotopy Appl. \textbf{27} (2025), no.~1, 107--122. \MR{4881650}

\bibitem[Hof98]{Hof98}
Detlev~W. Hoffmann, \emph{On the dimensions of anisotropic quadratic forms in
  {$I^4$}}, Inventiones Mathematicae \textbf{131} (1998), no.~1, 185--198.

\bibitem[Hou96]{Hou96}
Xiang-dong Hou, \emph{{${\rm GL}(m,2)$} acting on {$R(r,m)/R(r-1,m)$}},
  Discrete Math. \textbf{149} (1996), no.~1-3, 99--122. \MR{1375102}

\bibitem[HT98]{hoffmann199814}
Detlef~W Hoffmann and Jean-Pierre Tignol, \emph{On 14-dimensional quadratic
  forms in {$I^3$}, 8-dimensional forms in {$I^2$}, and the common value
  property}, Documenta Mathematica \textbf{3} (1998), 189.

\bibitem[IK00]{izhboldin2000some}
Oleg~T Izhboldin and Nikita~A Karpenko, \emph{Some new examples in the theory
  of quadratic forms}, Mathematische Zeitschrift \textbf{234} (2000), no.~4,
  647--695.

\bibitem[Kah05]{Kah05}
Bruno Kahn, \emph{On ``horizontal'' invariants attached to quadratic forms},
  Algebra and number theory, Hindustan Book Agency, Delhi, 2005, pp.~21--33.
  \MR{2193343}

\bibitem[Kar04]{Kar04}
Nikita~A. Karpenko, \emph{Holes in {$I^n$}}, Ann. Sci. \'{E}cole Norm. Sup. (4)
  \textbf{37} (2004), no.~6, 973--1002. \MR{2119244}

\bibitem[KLP68]{KasLinPet68}
Tadao Kasami, Shu Lin, and W.~Wesley Peterson, \emph{New generalizations of the
  {R}eed-{M}uller codes. {I}. {P}rimitive codes}, IEEE Trans. Inform. Theory
  \textbf{IT-14} (1968), 189--199. \MR{275989}

\bibitem[KRW72]{KneRosWar72}
Manfred Knebusch, Alex Rosenberg, and Roger Ware, \emph{Structure of {W}itt
  rings and quotients of abelian group rings}, American Journal of Mathematics
  \textbf{94} (1972), 119--155.

\bibitem[KT70]{KasTok70}
Tadao Kasami and Nobuki Tokura, \emph{On the weight structure of
  {R}eed--{M}uller codes}, IEEE Transactions on Information Theory \textbf{16}
  (1970), no.~6, 752--759.

\bibitem[KTA76]{KasTokAzu76}
Tadao Kasami, Nobuki Tokura, and Saburo Azumi, \emph{On the weight enumeration
  of weights less than 2.5d of {R}eed-{M}uller codes}, Information and Control
  \textbf{30} (1976), no.~4, 380--395. \MR{401324}

\bibitem[Lam05]{Lam05}
Tsit-Yuen Lam, \emph{Introduction to quadratic forms over fields}, Graduate
  Studies in Mathematics, vol.~67, American Mathematical Society, 2005.

\bibitem[Lor23]{Lor23}
Nico Lorenz, \emph{Pfister numbers over rigid fields}, Communications in
  Algebra \textbf{51} (2023), no.~3, 1158--1178.

\bibitem[Lor25]{Lor25}
\bysame, \emph{On generalized {A}lbert forms over discretely valued fields},
  Communications in Algebra \textbf{53} (2025), no.~12, 5194--5204.

\bibitem[McE69]{McE69}
Robert~J. McEliece, \emph{Quadratic forms over finite fields and second-order
  {R}eed--{M}uller codes}, JPL Space Programs Summary 37-58, Vol. III, Jet
  Propulsion Laboratory, 1969.

\bibitem[Mil70]{Mil70}
John Milnor, \emph{Algebraic {$K$}-theory and quadratic forms}, Invent. Math.
  \textbf{9} (1969/70), 318--344. \MR{260844}

\bibitem[MS77]{MacSlo77}
F.~J. MacWilliams and N.~J.~A. Sloane, \emph{The theory of error-correcting
  codes.}, North-Holland Mathematical Library, Vol. 16, North-Holland
  Publishing Co., Amsterdam-New York-Oxford, 1977. \MR{465510}

\bibitem[OVV07]{OVV}
D.~Orlov, A.~Vishik, and V.~Voevodsky, \emph{An exact sequence for
  {$K^M_\ast/2$} with applications to quadratic forms}, Ann. of Math. (2)
  \textbf{165} (2007), no.~1, 1--13. \MR{2276765}

\bibitem[Pel24]{Pel24}
Sarah Peluse, \emph{Finite field models in arithmetic combinatorics---twenty
  years on}, Surveys in combinatorics 2024, London Math. Soc. Lecture Note
  Ser., vol. 493, Cambridge Univ. Press, Cambridge, 2024, pp.~159--199.
  \MR{4749041}

\bibitem[PST09]{ParSurTig09}
R.~Parimala, V.~Suresh, and J.-P. Tignol, \emph{On the {P}fister number of
  quadratic forms}, Quadratic Forms---Algebra, Arithmetic, and Geometry,
  Contemporary Mathematics, vol. 493, American Mathematical Society,
  Providence, RI, 2009, pp.~327--338.

\bibitem[Rac13]{Rac13}
M{\'e}lanie Raczek, \emph{On the 3-{P}fister number of quadratic forms},
  Communications in Algebra \textbf{41} (2013), no.~1, 342--360.

\bibitem[Ros94]{Rost-14dim}
Markus Rost, \emph{On 14-dimensional quadratic forms, their spinors, and the
  difference of two octonion algebras}, Unpublished preprint, May 1994; revised
  September 1996 and March 1999. Available at
  \url{https://www.math.uni-bielefeld.de/~rost/14-dim.html}, 1994.

\bibitem[RS98]{Rains-Sloane98}
Eric~M. Rains and N.~J.~A. Sloane, \emph{Self-dual codes}, Handbook of coding
  theory, {V}ol. {I}, {II}, North-Holland, Amsterdam, 1998, pp.~177--294.
  \MR{1667939}

\bibitem[SB70]{SloBer70}
Neil J.~A. Sloane and Elwyn~R. Berlekamp, \emph{Weight enumerator for
  second-order {R}eed-{M}uller codes}, IEEE Trans. Inform. Theory
  \textbf{IT-16} (1970), 745--751. \MR{274196}

\bibitem[Scu24]{Scu24}
Stephen Scully, \emph{On the holes in {$I^n$} for symmetric bilinear forms in
  characteristic 2}, 2024, Preprint, \texttt{arXiv:2409.02061}.

\bibitem[Vis00]{Vis03}
A.~S. Vishik, \emph{On the dimensions of quadratic forms}, Dokl. Akad. Nauk
  \textbf{373} (2000), no.~4, 445--447. \MR{1788358}

\bibitem[Voe03]{Voe03}
Vladimir Voevodsky, \emph{Motivic cohomology with {$\Z/2$}-coefficients},
  Publications Mathematiques de l'IHES \textbf{98} (2003), 59--104.

\bibitem[Wad83]{Wad83}
Adrian~R. Wadsworth, \emph{$p$-{H}enselian fields: {$K$}-theory, {G}alois
  cohomology, and graded {W}itt rings}, Pacific Journal of Mathematics
  \textbf{105} (1983), no.~2, 473--496.

\end{thebibliography}
\end{document}